\documentclass[]{article}

\usepackage{amsmath,amssymb}
\usepackage{mathrsfs}
\usepackage{amsfonts}
\usepackage{leftidx}
\usepackage{enumerate}
\usepackage{amsthm}
\usepackage{scalerel}
\usepackage{MnSymbol}
\usepackage{hyperref}
\usepackage[authoryear]{natbib}

\title{Topological Vector Spaces: a non-standard approach with monads and galaxies.}
\author{Niels Charlier, Hans Vernaeve}

\theoremstyle{plain}
\newtheorem{thm}{Theorem}[section]
\newtheorem{lem}[thm]{Lemma}
\newtheorem{prop}[thm]{Proposition}
\newtheorem{cor}[thm]{Corrolary}
\newtheorem{rem}[thm]{Remark}
\theoremstyle{definition}
\newtheorem{defn}[thm]{Definition}
\newtheorem{nott}[thm]{Notation}
\newtheorem{exmp}[thm]{Example}

\newcommand{\compl}[1]{{#1}^\mathfrak{c}}

\newcommand{\Pwfin}{\mathcal{P}_\mathrm{fin}}
\newcommand{\Pw}{\mathcal{P}}

\newcommand{\Kappa}{\mathrm{K}}
\newcommand{\cldir}[1]{\vec{#1}}
\newcommand{\cldown}[1]{{}^\downarrow{#1}}
\newcommand{\cldownint}[1]{{}^{\downarrow}{#1}}
\newcommand{\cldownext}[1]{{}^{\mathrel{\rotatebox{-90}{$\rightsquigarrow$}}}{#1}}

\newcommand{\enlp}{\star}
\newcommand{\enl}{{}^\star}
\newcommand{\menl}{{}^\dagger}
\newcommand{\menlp}{\dagger}
\newcommand{\ees}{{}^\circ}
\newcommand{\ee}{{}^\ostar}
\newcommand{\se}{{}^\mathfrak{s}}
\newcommand{\sep}{\mathfrak{s}}
\newcommand{\ce}{{}^\mathfrak{c}}
\newcommand{\cep}{\mathfrak{c}}
\newcommand{\ps}{{}^{\filleddiamond}}
\newcommand{\stcp}{{}^\sigma}
\newcommand{\Bs}{\mathfrak{B}}
\newcommand{\R}{\mathbb{R}}
\newcommand{\C}{\mathbb{C}}
\newcommand{\K}{\mathbb{K}}
\newcommand{\N}{\mathbb{N}}

\newcommand{\U}{\mathcal{U}}
\newcommand{\V}{\mathcal{V}}

\newcommand{\B}{\mathcal{B}}
\newcommand{\Q}{\mathcal{Q}}
\newcommand{\I}{\mathcal{I}}
\newcommand{\m}{\mathbf{m}}

\newcommand{\nbh}{\mathrm{nbh}}
\newcommand{\fil}{\mathrm{fil~}}
\newcommand{\idl}{\mathrm{idl~}}
\newcommand{\fils}[1]{\mathrm{fil_{#1}~}}
\newcommand{\idls}[1]{\mathrm{idl_{#1}~}}
\newcommand{\cof}{\mathrm{cof}}
\newcommand{\coi}{\mathrm{coi}}
\newcommand{\mon}{\mathfrak{m}}
\newcommand{\gal}{\mathfrak{g}}
\newcommand{\iog}{\mathfrak{h}}

\newcommand{\M}{\mathfrak{M}}
\newcommand{\Ga}{\mathfrak{G}}
\newcommand{\mo}{\boldsymbol{\mu}}
\newcommand{\ga}{\boldsymbol{\Gamma}}

\newcommand{\f}{\mathfrak{f}}
\newcommand{\st}{\mathrm{st~}}

\newcommand{\stcom}{\hat{\mathrm{st}}~}
\newcommand{\stcomt}[1]{\hat{\mathrm{st}}_{#1}~}
\newcommand{\stt}[1]{\mathrm{st}_{#1}~}

\newcommand{\fin}{\mathrm{Fin}}
\newcommand{\gr}{\top}

\newcommand{\bdd}{\mathrm{Bdd}}
\newcommand{\comp}{\mathrm{comp}}

\newcommand{\pcomp}{\mathrm{pcomp}}

\newcommand{\ns}{\mathrm{ns}}

\newcommand{\pns}{\mathrm{pns}}
\newcommand{\bt}{\textnormal{b}}
\newcommand{\lpns}{\mathrm{lpns}}
\newcommand{\lpcomp}{\mathrm{lpcomp}}
\newcommand{\coco}{\mathrm{coco}}
\newcommand{\qs}{\mathrm{qs}}
\newcommand{\co}{\mathrm{co}}
\newcommand{\bal}{\mathrm{bal}}
\newcommand{\spann}{\mathrm{span}}
\newcommand{\cobal}{\mathrm{cobal}}

\newcommand{\accpoints}{\nabla} 
\def\restrict#1{\raise-.5ex\hbox{\ensuremath|}_{#1}}
\newcommand\restr[2]{{
		\left.\kern-\nulldelimiterspace 
		#1 
		\vphantom{\big|} 
		\right|_{#2} 
}}

\begin{document}

\maketitle

\begin{abstract}
A new and extensive formalism is developed for monads and galaxies in non-standard enlargements. It is shown that monads and galaxies can be manipulated using order-preserving and order-reversing set-to-set maps, and that set properties associated with these maps can be extended not only to internal sets but to all monads and galaxies. An abstract theory of Intersections of Galaxies is introduced. These concepts are applied to basic topology as well (locally convex) topological vector spaces, their various properties and completions, allowing these to be effortlessly defined and characterized. Duality theory is studied in this framework, allowing in particular to formulate new brief and insightful proofs for the theorems of Mackey-Arens and Grothendieck completeness without any technicalities.
\end{abstract}

\section*{Preliminaries}

Several proposed non-standard theories of topological vector spaces can be found in literature (\cite{stroyanlux,HensonMoore1972,YOUNG1972}). In this paper we introduce a new approach with new notations and more generalized definitions. The reader is only assumed to be familiar with the basic principles of non-standard analysis, in particular with saturated models.

\section{The general theory of monads and galaxies}

\subsection{Directedness and cofinality}

\begin{defn}
	A \textbf{preordered set} is a pair $(X, \leq)$ where $X$ is \textit{non-empty} set and $\leq$ a  reflexive, transitive relation on $J=X$. If the equivalence relation defined by $\leq \wedge \geq$ coincides with $=$, it is a \textbf{(partial) order}. If for any $a, b \in X: \exists c \in X: c \geq a, b$, the (pre-)order is \textbf{directed}. If for any $a,b \in J$, $a \leq b \vee a \geq b$, the (pre-)order is \textbf{total}. The \textbf{tails; heads; finite and infinite points} of $(X, \leq)$ are ($a \in [\enl]{X}$):
	\begin{displaymath}
	\begin{split}
	[\enl]{X}_{\leq a} := \lbrace x \in [\enl]{X}: x \leq a \rbrace~~~&~~~
	[\enl]{X}_{\geq a} := \lbrace x \in [\enl]{X}: x \geq a \rbrace\\
	\enl{X}_\f := \lbrace x \in \enl{X}: \exists a \in \stcp{X}: x \leq a \rbrace \supseteq \stcp{X};~~~&~~~
	\enl{X}_\infty := \lbrace \omega \in \enl{X}: \forall x \in \stcp{X}: \omega \geq x \rbrace.
	\end{split}
	\end{displaymath}
\end{defn}
If the greatest elements $X_\gr := \lbrace x \in {X}: \forall a \in X: x \geq a \rbrace = \emptyset$ then $\enl{X}_\infty \cap \enl{X}_\f = \emptyset$ and both sets are external. Otherwise $ \enl{X} = \enl{X}_\f \supseteq \enl{X}_\infty = \enl{X_\gr} $.

\begin{defn}
	Given preordered sets $(X, \leq)$ and $(Y, \leq)$ and a mapping $\psi: X \rightarrow Y$. Then $\psi$ is \textbf{order preserving}  if for any $x, x' \in X: y \geq y' \Rightarrow \psi(x) \geq \psi(x')$. It is \textbf{order reversing}  if for any $x, x' \in X: x \geq x' \Rightarrow \psi(x) \leq \psi(x')$.
\end{defn}

\begin{prop}
	Assume we are using a $\kappa^+$-saturated model. Given a pre-ordered set $(X, \leq)$ with $|X| \leq \kappa$. Then $X$ is directed iff $\enl{X}_\infty \neq \emptyset$.
	\begin{proof}
		Since $\enl{X}_\infty = \bigcap_{x \in X} \enl{X}_{\geq x}$, saturation implies that $\enl{X}_\infty \neq \emptyset$ iff $\lbrace \enl{X}_{\geq x} \rbrace_{x \in X}$ has non-empty finite intersections, which is equivalent to directedness. 
	\end{proof}
\end{prop}

\begin{prop} (Finite overspill) Assume we are using a $\kappa^+$-saturated model. Given a directed set $X$ with $|X| \leq \kappa$ and an internal $A \subseteq \enl{X}$ such that $\enl{X}_\infty \subseteq A$. Then there exists an $a \in X$ such that $\enl{X}_{\geq a} \subseteq A$.
	\begin{proof}
		Since $\bigcap_{x \in X} \enl{X}_{\geq x} \subseteq A$, i.e. $\bigcap_{x \in X} \enl{X}_{\geq x} \setminus A = \emptyset$, saturation implies $\lbrace \enl{X}_{\geq x} \setminus A \rbrace_{x \in X}$ has a non-empty finite intersection, which due to the directedness of $X$ implies the existence of a $x \in X$ such that $\enl{X}_{\geq x} \setminus A = \emptyset$, i.e. $\enl{X}_{\geq a} \subseteq A$.
	\end{proof}
\end{prop}

\begin{cor} (Infinite overspill) Assume we are using a $\kappa^+$-saturated model. Given a directed set $X$ with $|X| \leq \kappa$ and an internal $A \subseteq \enl{X}$ such that $\enl{X}_\f \subseteq A$. Then there exists an $\omega \in \enl{X}_\infty$ such that $\enl{X}_{\leq \omega} \subseteq A$.
	\begin{proof}
		Consider the internal set $Z = \lbrace \iota \in \enl{X}: \enl{X}_{\leq \iota} \not\subseteq A \rbrace$. Suppose $\enl{X}_\infty \subseteq Z$, then there exists a $a \in X$ such that $\enl{X}_{\geq a} \subseteq Z$, i.e. $\enl{X}_{\leq a} \not\subseteq A$, contradicting the assumption. Hence there exists an $\omega \in \enl{X}_\infty \setminus Z$, i.e. $\enl{X}_{\leq \omega} \subseteq A$.
	\end{proof}
\end{cor}

\begin{defn}
	Given a pre-ordered set $(X, \leq)$. Given two sets $A \subseteq \enl{X}$ and $B \subseteq \enl{X}$. Then $A$ is \textbf{cofinal} in $B$ iff for each $b \in B$ there exists a $a \in A$ such that $a \geq b$. $A$ is \textbf{coinitial} in $B$ iff for each $b \in B$ there exists $a \in A$ such that $a \leq b$. In the specific case where $A$ is also a subset of $B$ we denote this as $A \subseteq_\cof B$ for a cofinal subset and $A \subseteq_\coi B$ for a coinitial subset.
\end{defn}

\begin{defn}
	Given a pre-ordered set $(X, \leq)$. An \textbf{indexed family} in $X$ is a member of $X^J$ for some \textit{non-empty} index set $J$. Given a standard-indexed internal family $\mathscr{A} = (a_j)_{j \in \enl{J}}$, the standard points of $\mathscr{A}$ are 
\begin{displaymath}
\mathscr{A}_\sigma := \lbrace a_j: j \in \stcp{J} \rbrace.
\end{displaymath}
The pre-order of $X$ may be induced on $\mathscr{A}$. The notations for tails, heads, finite and infinite points are then generalized as follows:
\begin{displaymath}	
\begin{split}
\mathscr{A}_{\leq a} := \lbrace x \in \mathscr{A}: x \leq a \rbrace~~~&~~~
\mathscr{A}_{\geq a} := \lbrace x \in \mathscr{A}: x \geq a \rbrace\\
\mathscr{A}_\f := \lbrace x \in \mathscr{A}: \exists a \in \mathscr{A}_\sigma: x \leq a \rbrace \supseteq \mathscr{A}_\sigma;~~~&~~~
\mathscr{A}_\infty := \lbrace \omega \in \mathscr{A}: \forall x \in \mathscr{A}_\sigma: \omega \geq x \rbrace.
\end{split}
\end{displaymath}
\end{defn}

\begin{prop}
	Assume we are using a $\kappa^+$-saturated model. Given a pre-ordered set $(X,\leq)$ and a standard-indexed internal family $\mathscr{A} = (a_j)_{j \in \enl{J}}$ in $X$ with $|J| \leq \kappa$. Then $\mathscr{A}_\sigma$ is directed iff $\mathscr{A}_\infty \neq \emptyset$.
	
	\begin{proof}
		Consider the pre-order $\leq$ on $J$ induced by $\mathscr{A}$; i.e. $\forall j, j' \in J: j \leq j' \iff a_j \leq a_{j'}$. Then $\mathscr{A}_\sigma$ is directed iff $\stcp{J}$ is directed iff $\enl{J}_\infty \neq \emptyset$ iff $\mathscr{A}_{\infty} \neq \emptyset$ since $\mathscr{A}_{\infty} = \lbrace a_j \rbrace_{j \in \enl{J}_\infty}$.
	\end{proof}
\end{prop}

\begin{prop} \label{prop_cofincoinit_fam}
	Given a pre-ordered set $(X,\leq)$ and standard-indexed internal families $\mathscr{A} = (a_j)_{j \in \enl{J}}, \mathscr{B} = (b_k)_{k \in \enl{K}}$ in $X$, and internal $P \subseteq \lbrace a_j \rbrace_{j \in \enl{J}}, Q \subseteq \lbrace b_k \rbrace_{k \in \enl{K}}$. The following are equivalent:
	\begin{enumerate}[(i)]
		\item $\mathscr{A}_{\infty} \cap P$ coinitial in $\mathscr{B}_{\infty} \cap Q$.
		\item $\mathscr{A}_{\f} \cap P$ cofinal in $\mathscr{B}_{\f} \cap Q$ 		\big(if $P = \lbrace a_j \rbrace_{j \in \enl{J}}, Q = \lbrace b_k \rbrace_{k \in \enl{K}}$, this is equivalent to $\mathscr{A}_\sigma$ cofinal in $\mathscr{B}_\sigma$\big).
	\end{enumerate}	
	\begin{proof}	
		Consider the pre-order $\leq$ on $J$ (resp. $K$) induced by $\mathscr{A}$ (resp. $\mathscr{B}$); i.e. $\forall j, j' \in J: j \leq j' \iff a_j \leq a_{j'}$ (resp. $\forall k, k' \in K:k \leq k' \iff b_k \leq k_{k'}$). Note that $\mathscr{A}_{\infty} = \lbrace a_j \rbrace_{j \in \enl{J}_\infty}$ and $\mathscr{B}_{\infty} = \lbrace b_k \rbrace_{k \in \enl{K}_\infty}$.
					
		(i) $\Rightarrow$ (ii):  For any $q \in \mathscr{B}_\f \cap Q$, $\enl{J}_\infty \subseteq \lbrace j \in \enl{J}:\exists p \in P: q \leq p \leq a_j \rbrace$, so that because of overspill there exists a $j \in \enl{J}_\f$ and $p \in P$ such that $a_j \geq p \geq q$, which implies $p \in \mathscr{A}_{\f} \cap P$.  
		
		(ii) $\Rightarrow$ (i): Given $q \in \mathscr{B}_\infty \cap Q$, $\enl{J}_\f \subseteq \lbrace j \in \enl{J}: \exists p \in P: a_j \leq p \leq q \rbrace$, so that because of overspill there exists an $\omega \in \enl{J}_\infty$ and $p \in P$ such that $a_\omega \leq p \leq q$, which implies $p \in \mathscr{A}_{\infty} \cap P$.  . 
	\end{proof}
\end{prop}

\begin{cor}
	Given a pre-ordered set $(X,\leq)$ and directed $V, W \subseteq X$, the following are equivalent:
	\begin{enumerate}[(i)]
		\item $V$ cofinal in $W$.
		\item $\enl{V}_\infty$ coinitial in $\enl{W}_\infty$;		
	\end{enumerate}
	If $V \subseteq W$ this is equivalent to $\enl{V}_\infty \subseteq \enl{W}_\infty$ and $\enl{V}_\infty = \enl{V} \cap \enl{W}_\infty$.
	
	\begin{proof}
		For the last statement, notice that $\enl{V} \cap \enl{W}_\infty \subseteq \enl{V}_\infty$ always, hence $\enl{V}_\infty \subseteq \enl{W}_\infty$ implies $\enl{V}_\infty = \enl{V} \cap \enl{W}_\infty$. Now assume this to be the case and suppose there exists a $b \in W$ such that $\forall a \in V: a \not\geq b$. Because of the transfer principle $\enl{V} \cap \enl{X}_\infty = \emptyset$, a contradiction because $\enl{V}_\infty$ cannot be empty. Hence $V \subseteq_\cof W$ and thus $\enl{V}_\infty \subseteq_{[\coi]} \enl{W}_\infty$.
	\end{proof}
\end{cor}

\subsection{Monads and Galaxies}

From now on we work with a space called $E$. We will use indexed families - with cardinality limited by some $\kappa$ - of sets in $E$, accompanied by either one of two induced orders $\subseteq$ or $\supseteq$; as well as standard-indexed internal indexed families, denoted:
\begin{displaymath}
\begin{split}
\Bs_\kappa(E, \supseteq) := \bigcup_{|J| \leq \kappa, J \neq \emptyset} (\Pw(E), \supseteq)^J; ~~~&~~~ \ps\Bs_\kappa(E, \supseteq) := \bigcup_{|J| \leq \kappa, J \neq \emptyset} \enl(\Pw(E), \supseteq)^J;\\
\Bs_\kappa(E, \subseteq) := \bigcup_{|J| \leq \kappa, J \neq \emptyset} (\Pw(E), \subseteq)^J; ~~~&~~~ \ps\Bs_\kappa(E, \subseteq) := \bigcup_{|J| \leq \kappa, J \neq \emptyset} \enl(\Pw(E), \subseteq)^J.
\end{split}
\end{displaymath}
We may omit the subscript $\kappa$, in which case all cardinalities below saturation are included.

\begin{defn}
	Given $\U = (U_j)_{j \in \enl{J}} \in \ps\Bs(E, \supseteq)$, resp. $\ps\Bs(E, \subseteq)$. Then the indexed family $\cldir{\U}$ is the \textbf{$\cap$-closure}, resp. \textbf{$\cup$-closure} of $\U$, defined as
	\begin{displaymath}
		\cldir{\U} := (V_j)_{S \in \enl\Pwfin(J)}~~~V_{j_1, \dots, j_n} := \bigcap_{i = 1,\dots,n} U_{j_i}, \textnormal{resp.} \bigcup_{i = 1,\dots,n} U_{j_i}.
	\end{displaymath}
	Note that $\cldir{\U}_\sigma$ is always directed. For $\U \in {\ps\Bs(E, \supseteq)}$, $\bigcap\cldir\U_\sigma = \bigcap\U_\sigma$, $\emptyset \notin \cldir{\U}_\sigma$ iff $\bigcup \U_\infty \neq \emptyset$ iff $\U_\sigma$ has non-empty finite intersections. For $\U \in {\ps\Bs(E, \subseteq)}$, $\bigcup\cldir\U_\sigma = \bigcup\U_\sigma$, $\enl{E} \notin \cldir{\U}_\sigma$ iff $\bigcap \U_\infty \neq \enl{E}$ iff $\U_\sigma$ has non-covering finite unions.
\end{defn}

\begin{defn}
	Given $\U \in \Bs(E, \supseteq)$, resp. $\U \in \Bs(E, \subseteq)$. Then the \textbf{downward closure} $\cldown{\U}$ is the (unindexed) family
	\begin{displaymath}
	\cldown{\U} := \lbrace V \in \Pw(E): \exists U \in \U: V \supseteq U \textnormal{, resp. } V \subseteq U \rbrace.
	\end{displaymath}
	For $\U \in \ps\Bs(E, \supseteq)$, resp. $\U \in \ps\Bs(E, \subseteq)$, the \textbf{internal downward closure} $\cldownint\U$ is defined by transferring. The \textbf{external downward closure} $\cldownext{\U}$ is theng
	\begin{displaymath}
	\cldownext{\U} := \lbrace V \in \Pw(\enl{E}): U \in \U: V \supseteq U \textnormal{, resp. } V \subseteq U \rbrace.
	\end{displaymath}
\end{defn}

\begin{prop} (Cauchy principle) Given directed $\U \in \ps\Bs(E,\supseteq)$, resp. ${\ps\Bs(E,\subseteq)}$, then
	\begin{displaymath}
	\begin{split}
	\forall A \in \enl{\Pw}(E):~~~& \bigcup \U_\infty \subseteq A \iff \exists U \in \U_\sigma: U \subseteq A;  \\	resp.~~&A \subseteq \bigcap \U_\infty \iff \exists U \in \U_\sigma: A \subseteq U.
	\end{split}
	\end{displaymath}
	\begin{proof}
		This follows from finite overspill.
	\end{proof}
\end{prop}

\begin{prop} (Nucleus principle) Given directed $\U \in \ps\Bs(E,\supseteq)$, resp. $\ps\Bs(E,\subseteq)$, then
	\begin{displaymath}
	\begin{split}
	\forall A \in \enl{\Pw}(E):~~~& A \subseteq \bigcap \U_\sigma \iff \exists U \in \U_\infty: A \subseteq U;  \\	resp.~~&\bigcup\U_\sigma \subseteq A  \iff \exists U \in \U_\infty:  U \subseteq A.
	\end{split}
	\end{displaymath}
	\begin{proof}
		This follows from infinite overspill.
	\end{proof}
\end{prop}

\begin{prop} \label{prop_orderbasis}	
	Given directed $\U \in \ps\Bs(E,\supseteq)$, resp. ${\ps\Bs(E,\subseteq)}$. The following statements are equivalent:
	
	\begin{enumerate}[(i)]
		\item $\bigcap \U_\sigma \subseteq \bigcap \V_\sigma$, resp. $\bigcup \U_\sigma \supseteq \bigcup \V_\sigma$;
		\item $\bigcup \U_\infty \subseteq \bigcup \V_\infty$, resp. $\bigcap \U_\infty \supseteq \bigcap \V_\infty$;
		\item $\U_\sigma$ cofinal in $\V_\sigma$;
		\item $\V_\infty$ coinitial in $\U_\infty$.
	\end{enumerate}
	
	\begin{proof}
		The equivalence (iii) $\iff$ (iv) follows from proposition \ref{prop_cofincoinit_fam}. The implications (iii) $\Rightarrow$ (i) and (iv) $\Rightarrow$ (ii) are clear. The implications (ii) $\Rightarrow$ (iii) and (i) $\Rightarrow$ (iv) follow from the Cauchy, resp. nucleus principle.
	\end{proof}
\end{prop}

\begin{prop} \label{prop_ordermonad} Given directed $\U \in \ps\Bs(E,\supseteq)$, resp. ${\ps\Bs(E,\subseteq)}$. 
	\begin{displaymath}	
		\bigcap \U_\sigma = \bigcap \U_\f = \bigcup \U_\infty
 \textnormal{, resp. }
		\bigcup \U_\sigma = \bigcup \U_f = \bigcap \U_\infty.
	\end{displaymath}
	
	\begin{proof}
		The inclusion $\bigcap \U_\sigma \subseteq \bigcup \U_\infty$ follows from the nucleus principle; $\bigcap \U_\infty \subseteq \bigcup \U_\sigma$ from the Cauchy principe (each time applied to a singleton). The other inclusions are clear.
	\end{proof}	
\end{prop}

\begin{defn}
	Proposition \ref{prop_orderbasis} defines a pre-order $\succeq$ on the directed families of $\ps\Bs(E,\supseteq)$, resp. $\ps\Bs(E,\subseteq)$. This pre-order is extended to all families by letting $\U \succeq \V \iff \cldir{\U} \succeq \cldir{\V}.$ We will denote $\sim$ for the equivalence relationship on $\ps\Bs(E,\subseteq)$ defined by $\preceq \land \succeq$.
	
\end{defn}

\begin{defn}
	
	The \textbf{monads} and \textbf{galaxies} [of cofinality $\leq \kappa$] in $E$ are
	\begin{displaymath}
	\begin{split}
	\M_{[\kappa]}(\enl{E}) := \Bs_{[\kappa]}(E, \supseteq) / {\sim};~~~&~~~~
	\ps\M_{[\kappa]}(\enl{E}) := \ps\Bs_{[\kappa]}(E, \supseteq) / {\sim};\\
	\Ga_{[\kappa]}(\enl{E}) := \Bs_{[\kappa]}(E, \subseteq) / {\sim};~~~&~~~~
\ps\Ga_{[\kappa]}(\enl{E}) := \ps\Bs_{[\kappa]}(E, \subseteq) / {\sim};
	\end{split}
	\end{displaymath}	
\end{defn}

We identify the left-side, the \textit{(standard)} monads and galaxies, as subsets of the right-side, the $\ps$monads and $\ps$galaxies. A representant of a $[\ps]$monad or $[\ps]$galaxy is called a $[\ps]$\textbf{(sub)base} - a subbase $\U$ is a base if $\U_\sigma$ is directed. The pre-order $\succeq$ induces an order on the quotients, where the left (resp. right) side of the order relation is called a \textbf{finer} (resp. \textbf{coarser}) $[\ps]$monad/galaxy. Consider the mappings - the representant $\U$ is a base (using proposition \ref{prop_ordermonad}):
\begin{displaymath}
\begin{split}
\ps\M({\enl{E}}) \rightarrow \Pw(\enl{E}):  \U \rightarrow \bigcup \U_\infty = \bigcap \U_\sigma;\\
\ps\Ga(\enl{E}) \rightarrow \Pw(\enl{E}):  \U \rightarrow \bigcap \U_\infty = \bigcup \U_\sigma.
\end{split}
\end{displaymath}  
Because of proposition \ref{prop_orderbasis} these are in fact injections, and we will consistently use them as implicit identifications. In this identification the `is finer than' order corresponds with $\subseteq$ for monads and $\supseteq$ for galaxies.
 Then we have:
\begin{displaymath}
\begin{gathered}
\M_1(\enl{E}) (= \M_2(\enl{E}) = \dots) = \Ga_1(\enl{E}) (= \Ga_2(\enl{E}) = \dots) = \stcp{\Pw}(E);\\
\ps\M_1(\enl{E}) (= \ps\M_2(\enl{E}) = \dots) = \ps\Ga_1(\enl{E}) (= \ps\Ga_2(\enl{E}) = \dots) = \enl{\Pw}(E).
\end{gathered}
\end{displaymath}
However, all other monads and galaxies are identified as external sets. We must thus be careful: $\M(\enl{E})$ and $\Ga(\enl{E})$ are standard sets, but the implicit identification of their elements as sets is external. 

\begin{nott}
	Given a set $A \subseteq \enl{E}$ we denote $\M(A)$, resp. $\Ga(A)$ for all monads, resp. galaxies contained in $A$.
\end{nott}

\begin{prop}
	$\M(\enl{E}) \cap \ps\M_{\kappa}(\enl{E}) = \M_{\kappa}(\enl{E})$; $\Ga(\enl{E}) \cap \ps\Ga_{\kappa}(\enl{E}) = \Ga_{\kappa}(\enl{E})$.
	
	\begin{proof}
		Given $\mon \in \M(\enl{E}) \cap \ps\M_{\kappa}(\enl{E})$.
		 Suppose $\U$ is a standard base for $\mon$ and $(U_j)_{j \in \enl{J}}$ is a $\ps$base with $|J| \leq \kappa$. Then for any $j \in \stcp{J}$, let $V_j$ be a set in $\U$ such that $\enl{V_j} \subseteq U_j$. Then $(V_j)_{j \in J}$ is a standard base for $\mon$. For galaxies the proof is analogue.
	\end{proof}
\end{prop}

\begin{exmp}
	Standard copies of sets are galaxies with finite sets as a base. In particular, $\stcp{E}$ is the galaxy that has all finite sets of $E$ as its base.
	For any pre-ordered set $X$, $\enl{X}_\f$ is a galaxy and $\enl{X}_\infty$ is a monad.
\end{exmp}

\begin{defn}
	We call $\emptyset$, resp. $\enl{E}$ the \textbf{trivial} monad, resp. galaxy.	A $\ps$monad or $\ps$galaxy is \textbf{proper} if it is not trivial. The cauchy principle implies that a monad., resp. galaxy with a given (sub)base is proper iff the (sub)base has non-empty finite intersections, resp. non-covering finite unions.
\end{defn}

\begin{defn}
Given spaces $E_1, \dots, E_n, F$ and an internal mapping $
\psi: \enl\mathcal{P}(E_1) \times \dots \times \enl\mathcal{P}(E_n) \rightarrow \enl\mathcal{P}(F)
$
 that is either order preserving in each argument or order reversing in each argument. We define, given $\ps$bases $(U_{j_1})_{j_1 \in \enl{J_1}}$ in $E_1, \dots, (U_{j_n})_{j_n \in \enl{J_n}}$ in $E_n$:
\begin{displaymath}
\begin{gathered}
\menl\psi((U_{j_1})_{j_1 \in \enl{J_1}}, \dots, (U_{j_n})_{j_n \in \enl{J_n}} ) := (\psi(U_{j_1}, \dots U_{j_n}))_{(j_1, \dots, j_n) \in \enl{(J_1 \times \dots \times J_n)}}.
\end{gathered}
\end{displaymath}
Now note that order preservation or reversion and proposition \ref{prop_orderbasis} imply that given $\ps$bases $\U_1, \dots \U_n, \V_1, \dots \V_n$ we have $\U_1 \sim V_1,\dots,\U_n \sim \V_n \Rightarrow \menl\psi(\U_1, \dots, \U_n) \sim \menl\psi(\V_1, \dots, \V_n$). Hence we may interpret the definition above of $\menl\psi$ as a mapping of monads and galaxies that is independent of the representant. If the mapping is order preserving, $\prod_{j \leq n} \ps{\M}_\kappa(E_j)$ is mapped to $\ps{\M}_\kappa(E)$ and $\prod_{j \leq n} \ps{\Ga}_\kappa(E_j)$ is mapped to $\ps{\Ga}_\kappa(E)$. If the mapping is order reversing, $\prod_{j \leq n} \ps{\M}_\kappa(E_j)$ is mapped to $\ps{\Ga}_\kappa(E)$ and $\prod_{j \leq n} \ps{\Ga}_\kappa(E_j)$ is mapped to $\ps{\M}_\kappa(E)$. If $\psi$ is a standard map, then standard monads and galaxies are mapped onto standard monads and galaxies.
\end{defn}

\begin{prop} \label{prop_imagemonadset}
	Given spaces $E_1, \dots, E_n, F$ and a mapping $\psi: \Pw(\enl{E_1}) \times \dots \times \mathcal{P}(\enl{E_n}) \rightarrow \Pw(\enl{F})$ that is order preserving or reversing, such that $\psi$ restricted to $\enl(\Pw(E_1) \times \dots \times \Pw(E_n))$ is internal. For any $\mon_1 \in \ps{\M}(E_1),\dots,\mon_n \in \ps{\M}(E_n)$ or $\gal_1 \in \ps{\Ga}(E_1),\dots,\gal_n \in \ps{\Ga}(E_n)$, we have (applying $\menlp$ to the internal restriction):
	\begin{displaymath}
		\psi(\mon_1,\dots\mon_n) = \menl\psi(\mon_1, \dots\mon_n) ~~~~~~ \psi(\gal_1, \dots, \gal_n) = \menl\psi(\gal_1, \dots, \gal_n).
	\end{displaymath}
	
	\begin{proof}
		Assume a single argument mapping. Given monad $\mon$ with $\ps$base $\U$ (the proof for a galaxy is analogue). If $\psi$ is order-preserving we have 
		\begin{displaymath}
		\begin{split}
		\menl\psi(\mon) &= \bigcap \psi(\U_\sigma) \supseteq \psi(\bigcap \U_\sigma) = \psi(\mon) = \psi(\bigcup \U_\infty) \supseteq \bigcup \psi(\U_\infty) = \menl\psi(\mon).		
		\end{split}
		\end{displaymath}
		while if $\psi$ is order-reversing 
		\begin{displaymath}
		\begin{split}
		\menl\psi(\mon) &= \bigcup \psi(\U_\sigma) \supseteq \psi(\bigcup \U_\sigma) = \psi(\mon) = \psi(\bigcup \U_\infty) \supseteq \bigcap \psi(\U_\infty) = \menl\psi(\mon).		
		\end{split}
		\end{displaymath}
		For multiple arguments, first fixate all but one argument on an arbitrary internal set; then use induction for multiple monads or galaxies.
	\end{proof}
\end{prop}

\begin{cor} \label{cor_imagemonad}
	Given spaces $E_1, \dots, E_n, F$ and an internal mapping $\psi: \enl{E_1} \times \dots \times \enl{E_n} \rightarrow \enl{F}$. For any $\mon_1 \in \ps{\M}(E), \dots, \mon_n \in \ps{\M}(E_n)$ or $\gal_1 \in \ps{\Ga}(E_1), \dots, \gal_n \in \ps{\Ga}(E_n)$ we have
	\begin{displaymath}
	\psi(\mon_1,\dots\mon_n) = \menl\psi(\mon_1, \dots\mon_n) ~~~~~~ \psi(\gal_1, \dots, \gal_n) = \menl\psi(\gal_1, \dots, \gal_n)
	\end{displaymath}
	where we have interpreted $\psi$ as a $\mathcal{P}(\enl{E_1}) \times \dots \times \mathcal{P}(\enl{E_n}) \rightarrow \mathcal{P}(\enl{G})$ mapping that is internal for internal sets and order-preserving. 
\end{cor}

\begin{nott}
	For $V \subseteq E$, respectively $A \subseteq \enl{E}$ we denote $\compl{V} := E \setminus V$, respectively $\compl{A} := \enl{E} \setminus A$.
	
\end{nott}

\begin{cor}
	$\ps{\M}(E)$ and $\ps{\Ga}(E)$ are each closed for finite union and intersection.\footnote{In fact, one can prove that $\M(\enl{E})$, resp. $\Ga(\enl{E})$ are closed for infinite intersections, resp. unions; but this is not the case for unions of monads and intersections of galaxies.} The complement of a $\ps$monad, respectively $\ps$galaxy is a $\ps$galaxy, respectively $\ps$monad. 	
	\begin{proof}
		Consider the mappings ($A, B \subseteq \enl{E}$): $(A, B) \rightarrow A \cap B$, $(A, B) \rightarrow A \cup B$ and $A \rightarrow \compl{A}$. The first two are order preserving in each argument; the latter is order reversing and each is internal for internal sets; hence their $\menlp$-extensions coincides with the general set operators.
	\end{proof}
\end{cor}

\begin{defn}
	Given an internal order-preserving mapping $\psi: \enl\mathcal{P}(E) \rightarrow \enl\mathcal{P}(E)$. Then the set property as defined by the set
	\begin{displaymath}
	\Q_\psi = \lbrace V \in \enl\Pw(E): \psi(V) = V \rbrace.
	\end{displaymath}
	is extendable to a \textbf{$\menlp$-property} as follows:
	\begin{displaymath}
	\menl\Q_\psi := \lbrace \mathfrak{a} \in \ps{\M}(E) \cup \ps{\Ga}(E): \menl\psi(\mathfrak{a}) = \mathfrak{a} \rbrace.
	\end{displaymath}
	This is equivalent to $\mathfrak{a}$ having a $\ps$base in $\Q_\psi$. If $\psi$ is defined on external sets, then the property is also defined on external sets by 
	\begin{displaymath}
	\ees{\Q}_\psi := \lbrace V \in \Pw(\enl{E}): \psi(V) = V \rbrace.	
	\end{displaymath}
	In this case $\menl\Q_\psi$ is precisely the restriction of this set to monads and galaxies and we may omit the $\menlp$ symbol.
\end{defn}

\begin{lem} \label{lem_hullintersect}
	Given an order-preserving map $\psi: \Pw(\enl{E}) \rightarrow \Pw(\enl{E})$ such that $\psi(A) \supseteq A$ (respectively $\psi(A) \subseteq A$) for any $A \subseteq \enl{E}$ and associated set property $\ees{\Q}_\psi$. Then $\ees{\Q}_\psi$ is closed for intersection (respectively union).
	
	\begin{proof}
		Given a family of sets $\lbrace U_j \rbrace_{j \in J} \subseteq \ees{\Q}_\psi$ and $U$ its intersection. Then
		\begin{displaymath}
		\psi(U) = \psi \big(\bigcap_{j \in J} U_j\big)
		\subseteq \bigcap_{j \in J} \psi(U_j)
		= \bigcap_{j \in J} U_j = U.
		\end{displaymath}
		If the assumption was $\psi(U) \supseteq U$, the result follows. The other case is analogue.
	\end{proof}
\end{lem}

\begin{rem}
	The order-preserving extension of an internal mapping is unique on monads and galaxies, however this is not the case for other external sets. In fact, we will later discover an example where two order-preserving mappings coincide entirely on standard sets, internal sets, monads and galaxies, but not on all external sets.
\end{rem}

\begin{prop} (Combined Cauchy Principle) Given $\mon \in \ps\M(\enl{E})$ and $\gal \in \ps\Ga(\enl{E})$. Then
	\begin{displaymath}
	\mon \subseteq \gal \iff \exists A \subseteq \enl{\mathcal{P}(E)}: \mon \subseteq A \subseteq \gal. 
	\end{displaymath}
	\begin{proof}
		The inclusion is equivalent to $\mon \cap \compl{\gal}$ being a trivial monad; i.e. each $\ps$base of this monad containing the empty set. Given $\ps$bases $\U$ for $\mon$ and $\V$ for $\gal$, the statement is also equivalent with the existence of a $U \in \U$ and $V \in \V$ such that $U \subseteq V$.
	\end{proof}
\end{prop}

\begin{nott}
	Given a $\enl$net $\mathscr{A} = (a_j)_{j \in \enl{J}}$, i.e. $(J, \leq)$ is a directed pre-ordered standard set. Then we denote ($\iota \in J$):
	\begin{displaymath}
	\begin{gathered}	
	\mathscr{A}_{\enlp} :=  \lbrace a_j \rbrace_{j \in \enl{J}};~~~~~
	\mathscr{A}_{\leq \iota} := \lbrace a_j \rbrace_{j \leq \iota};~~~~~
	\mathscr{A}_{\geq \iota} := \lbrace a_j \rbrace_{j \geq \iota};\\
	\mathscr{A}_\sigma := \lbrace a_j \rbrace_{j \in \stcp{J}};~~~~~
	\mathscr{A}_\f := \lbrace a_j \rbrace_{j \in \enl{J}_\f};~~~~~~
	\mathscr{A}_\infty :=  \lbrace a_j \rbrace_{j \in \enl{J}_\infty}.
	\end{gathered}
	\end{displaymath}
\end{nott}

Given $\mathscr{A} = (a_j)_{j \in \enl{J}}$ with $|J| \leq \kappa$. Then $(\mathscr{A}_{\geq j})_{j \in \enl{J}}$ is a $\ps$base (of tails), hence $\mathscr{A}_\infty \in \ps\M_\kappa(\enl{E})$. 

\begin{defn}
	The \textbf{$\kappa$-thin} $[\ps]$monads are the infinite point sets of $[\enl]$nets with a standard index of cardinality $\leq \kappa$.  We denote them as $[\ps]\M_{.\kappa}(\enl{E}) \subseteq [\ps]\M_\kappa(\enl{E})$. In particular, we call $\aleph_0$-thin $[\ps]$monads \textbf{sequential}, they are associated with $[\enl]$sequences.
\end{defn}

The following propositions show that nets may be used as representatives of monads and sequences as representatives of monads with countable cofinality.

\begin{prop}
	Given $\mon \in \M_\kappa(\enl{E})$, there exists a cardinality $\lambda \geq \kappa$ such that $\mon \in \M_{.\lambda}(\enl{E})$.
	
	\begin{proof} Let $\U$ be a base for $\mon$.
		Let $J = \lbrace (U, e): U \in \U, e \in U \rbrace$ with $(U, e) \leq (U', e') \iff U \leq U'$ and consider the net $\mathscr{A} = (a_j)_{j \in J}; a_{(U, e)} := e$. Then $\enl\mathscr{A}_\infty = \mon$.
	\end{proof}
\end{prop}

\begin{prop}
	Given $\mon \in [\ps]\M_\kappa(\enl{E})$, there exists an $\mon' \in [\ps]\M_{.\kappa}(\mon)$.
	\begin{proof}		
		Consider the standard case: if $\U = (U_j)_{j \in J}$ is a base for $\mon$, choose $a_j \in U_j$ for each $j \in J$ and let $\mathscr{A} = (a_j)_{j \in J}$; then $\enl\mathscr{A}_\infty \subseteq \mon$. Then because of the transfer principle we can do the same for a $\ps$base.\footnote{Note that we are applying the \textit{transferred} axiom of choice rather than the axiom of choice itself on $\enl{J}$ to ensure the choice is internal.}
	\end{proof}
\end{prop}

\begin{cor} \label{cor_countablemonadseq}
	Any monad of countable cofinality contains a sequential finer monad.
\end{cor}

\subsection{Properties of standard monads and galaxies}

\begin{defn}
	A \textbf{ring of sets} $\Sigma$ over $E$ is a family of sets $\Sigma$ in $E$ such that $\emptyset \in \Sigma, E \in \Sigma$, closed for finite intersection and union. The \textbf{discrete ring} is the ring of all subsets of $E$.
\end{defn}

\begin{defn} 
	Given a ring $\Sigma$ and a set $A \in \enl{E}$. Then the \textbf{filter} and \textbf{monad} over $\Sigma$ generated by $A$ are
	\begin{displaymath}
	\begin{split}
		\fils{\Sigma}{A} &:= \lbrace U \in \Sigma, A \subseteq \enl{U} \rbrace;\\
		\mo_\Sigma(A) &:= \bigcap \stcp(\fils{\Sigma}{A}) = \bigcup \enl(\fils{\Sigma}{A})_\infty.	
	\end{split}
	\end{displaymath} 
	The \textbf{ideal} and \textbf{galaxy} over $\Sigma$ generated by $A$ are
	\begin{displaymath}
	\begin{split}
	\idls{\Sigma}{A} &:= \lbrace U \in \Sigma, \enl{U} \subseteq A \rbrace;\\
	\ga_\Sigma(A) &:= \bigcup \stcp(\fils{\Sigma}{A}) = \bigcap \enl(\fils{\Sigma}{A})_\infty.	
	\end{split}
	\end{displaymath} 
	
	In the absence of an explicit $\Sigma$, the discrete ring is implied. In the case of a singleton, we may omit the braces.	
\end{defn}

Now notice that for $\mon \in \M(\enl{E})$ with base $\U$, we have $\fil \mon = \cldown\U$ and $\mo(\mon) = \mon$, hence $\mon \rightarrow \fil \mon$ is a standard bijection into the set of filters on $E$. The filter is the greatest base (for the $\subseteq$ order) of $\mon$; the bases of $\mon$ are precisely the cofinal subsets of its filter. All properties of (standard) monads can be interpreted as properties of filters, relating them to standard analysis. The case for galaxies and ideals is analogue.

\begin{defn}
	A non-trivial monad that does not contain any strictly finer non-trivial monads is an \textbf{ultramonad}. A non-trivial galaxy that is not contained in any strictly finer non-trivial galaxies is an \textbf{ultragalaxy}.
\end{defn}

\begin{prop}
	The ultramonads of $E$ are precisely $\lbrace \mo(e): e \in \enl{E} \rbrace$.
	\begin{proof}
		We may construct a base for a strictly finer non-trivial monad iff there exists a $U \subseteq E$ such that $\fil \mon$ does not contain either $U$ or its complement.
		First consider $\mon = \mo(e)$ for $e \in \enl{E}$, then $U \in \fil \mon$ iff $e \in \enl{U}$ but otherwise $e \in \enl{\compl{U}}$. For the opposite implication, suppose $\mon$ contains points $e, e'$ such that $e' \notin \mo(e)$, then $\mo(e)$ is a strictly finer non-trivial monad of $\mon$.
	\end{proof}
\end{prop}

Similarly, the ultragalaxies are precisely $\lbrace \ga(\compl{e}) = \compl{\mo(e)}: e \in \enl{E} \rbrace$.

\begin{defn}
	The \textbf{kernel} of a monad $\mon \in \M(\enl{E})$ is
	\begin{displaymath}
	\ker \mon := \bigcap_{U \in \fil{\mon}} U = \mon \cap \stcp{E}.
	\end{displaymath}	
	Note that $\enl\ker \mon = \bigcap \enl(\fil \mon)_\infty \subseteq \mon$. A monad is \textbf{free} if $\ker \mon = \emptyset$. A monad is \textbf{principal} if $\enl\ker \mon = \mon$ (which is equivalent to finite cofinality).  An ultramonad $\mo(e)$ is principal iff $e$ is standard and free otherwise.
\end{defn}

\begin{defn}
	A galaxy $\gal \in \M(\enl{E})$ is \textbf{covering} if $\stcp{E} \subseteq \gal$.
\end{defn}

Hence free monads and covering galaxies are each others complements.

\begin{defn} Given a ring $\Sigma$. A set $A \subseteq \enl{E}$ is \textbf{$\Sigma$-chromatic} if $\forall a \in A: \mo_\Sigma(a) \subseteq A$. Naturally, for any $U \in \Sigma$, $\enl{U}$ is $\Sigma$-chromatic. Again the absence of an explicit $\Sigma$ refers to the discrete ring so that a set is \textbf{chromatic} if it is a union of monads.
\end{defn}

The following is a generalization of a known result (\cite{stroyanlux}, 8.2.2, p. 199)\footnote{Furthermore, this proof does not use the axiom of choice.}:

\begin{prop} 
	Given a ring $\Sigma$ and $\mon \in \ps\M(\enl{E})$. Then
		 \begin{displaymath}
		 \mo_\Sigma(\mon) = \bigcup_{e \in \mon} \mo_\Sigma(e)
		\end{displaymath}
	\begin{proof}
		The inclusion $\bigcup_{e \in \mon} \mo_\Sigma(e) \subseteq \mo_\Sigma(\mon)$ is clear, assume then that it is strict.	Take $e \in \mo_\Sigma(\mon) \setminus \mon$. Let $\mathfrak{A} := \lbrace \mo_\Sigma(e') : e' \in \mon \rbrace \subseteq \M(\enl{E})$ and note that for any $\mathfrak{a} \in \mathfrak{A}$, $\mathfrak{a} \subseteq \compl{\lbrace e \rbrace}$. Using the Cauchy principle, consider the galaxy $\gal_e$ with base $\U = \lbrace U \in \Sigma: \exists \mathfrak{a} \in \mathfrak{A}: \mathfrak{a} \subseteq \enl{U}; e \notin \enl{U} \rbrace$. Since we have $\mon \subseteq \gal_e$, the combined Cauchy principle implies that $\mon \subseteq \enl{U}$ with $U \in \U \subseteq \Sigma$ and $e \not\in \enl{U}$. This is a contradiction because $e \in \mo_\Sigma(\mon)$. Hence we can conclude that $\mon = \mo_\Sigma(\mon)$. 
		
	\end{proof}
\end{prop}

\begin{cor} \label{cor_mongal_chrom}
	Any $\mon \in \ps\M(\enl{E})$ or $\gal \in \ps\Ga(\enl{E})$ is $\Sigma$-chromatic iff it is standard and has a base in $\Sigma$.
	
	\begin{proof}
		If $\mon$ is standard with a base in $\Sigma$, then clearly $\mon = \mo_\Sigma(\mon)$. If otherwise $\mon$ is $\Sigma$-chromatic then $\mon = \bigcup_{e \in \mon} \mo_\Sigma(e) = \mo_\Sigma(\mon)$. 
		If $\gal$ is a standard galaxy with a base in $\Sigma$, it is a union of $\Sigma$-chromatic sets and therefore $\Sigma$-chromatic. Assume otherwise that $\gal$ is $\Sigma$-chromatic. Since the complement of a chromatic set is chromatic, the first part of the proof implies that $\gal$ is standard. Take any $U \in \idl \gal$, then $\mo_\Sigma(\enl{U}) = \bigcup_{e \in \enl{U}} \mo_\Sigma(e) \subseteq \gal$; hence because of the combined Cauchy principle, there exists a $V \in \Sigma \cap \idl \gal$ such that $U \subseteq V$.
	\end{proof}
\end{cor}

\begin{cor} \label{cor_internalchromaticstandard}
	Given a ring $\Sigma$ and an internal set $A$. Then $A$ is $\Sigma$-chromatic iff $A \in \stcp{\Sigma}$. In particular, $A$ is chromatic iff it is standard.
\end{cor}

\begin{defn}
	Given a ring $\Sigma$ and $A, B \subseteq \enl{E}$. Then $A$ is \textbf{$\Sigma$-distinguishable} from $B$ iff $B \cap \mo_\Sigma(A) = \emptyset$ (equivalently due to the Cauchy principle, there exists a $V \in \Sigma$ such that $A \subseteq \enl{V}$ and $B \subseteq \enl{\compl{V}}$). Then $A$ and $B$ are \textbf{$\Sigma$-separated} iff $\mo_\Sigma(A) \cap \mo_\Sigma(B) = \emptyset$ (equivalently due to the Cauchy principle, there exists disjunctive $V, W \in \Sigma$ with met $A \subseteq V$ and $B \subseteq W$). We apply these terms also to points $e, e' \in \enl{E}$ (as singletons). Note that a pair being mutually distinguishable is strictly weaker than being separated.
\end{defn}

\begin{nott}
	Given a ring $\Sigma$, then the complements of sets in $\Sigma$ are a ring as well, denoted as $\overline{\Sigma}$. The discrete ring is of course equal to its complement.
\end{nott}

\begin{lem} \label{lem_equivdis}
	Given a ring $\Sigma$, and $e, e' \in \enl{E}$:
	$e \in \mo_\Sigma(e') \iff e' \in \mo_{\overline\Sigma}(e)$.
	\begin{proof}
		Suppose $e \notin \mo_{\overline\Sigma}(e')$, then there exists a $U \in \overline\Sigma$ with $e \in U$ and $e' \notin U$. Then $\compl{U} \in \Sigma$, hence $e' \notin \mo_\Sigma(e')$. The opposite implication is identical.
	\end{proof}
\end{lem}

\begin{lem} \label{lem_equivcomp}
	Given a ring $\Sigma$. Then $A \subseteq \enl{E}$ is $\Sigma$-chromatic iff $\compl{A}$ is $\overline{\Sigma}$-chromatic.
	\begin{proof}		
		Given $e \in \compl{A}$, then  $e \notin \mo_\Sigma(e')$, then because of lemma \ref{lem_equivdis}, $e' \notin \mo_{\overline\Sigma}(e)$ for any $e' \in A$, hence $\mo_{\overline\Sigma}(e) \subseteq \compl{A}$. The opposite implication is identical.
	\end{proof}
\end{lem}

\begin{prop} \label{prop_chromgal}
	Given a ring $\Sigma$. Then $A \subseteq \enl{E}$ is $\Sigma$-chromatic iff for all $e \in \compl{A}: A \subseteq \ga_{\Sigma}(\compl{e})$.	
	\begin{proof}
		Because of lemma \ref{lem_equivcomp}, $A \subseteq \enl{E}$ is \textbf{$\Sigma$-chromatic} iff
		\begin{displaymath}		
		\compl{A} = \bigcup_{e \in A} \mo_{\overline{\Sigma}}(e) \iff 
		A = \compl{\big( \bigcup_{e \in A} \mo_{\overline{\Sigma}}(e) \big)}
		= \bigcap \compl{\mo_{\overline{\Sigma}}(e)} = \bigcap \ga_{\Sigma}(\compl{e}).
		\end{displaymath}
	\end{proof}
\end{prop}

In particular, a set is chromatic iff it is exactly equal to the union of all standard monads it contains iff it it is exactly equal to the intersection of all standard galaxies that contain it.

\begin{lem} \label{lem_chrominclusion} (Chromatic inclusion)
	Given two sets $A, B \subseteq \enl{E}$ where $A$ is chromatic. Then	
	\begin{displaymath}
		A \subseteq B \iff \M(A) \subseteq \M(B).
	\end{displaymath}
	\begin{proof}
		The $\Rightarrow$ implication is trivial. For $\Leftarrow$, take $e \in A$. Since $A$ is chromatic $\mo(e) \subseteq A$. Then $\mo(e) \in \M(A) \subseteq \M(B)$, hence, $\mo(e) \subseteq B$.
	\end{proof}
\end{lem}

\begin{prop} \label{prop_sequentialmonads}
	Given $\mon \in \M_{.\aleph_0}(\enl{E})$, $Q \subseteq \enl{E}$ such that $\stcp{E} \subseteq Q$. Then $\mon \subseteq Q \iff \mon \subseteq \ga(Q)$.
	
	\begin{proof}
		Let $\mathscr{A} = \enl(a_n)_{n \in \N}$ be a sequence such that $\mon = \mathscr{A}_\infty$. Assume $\mon \subseteq Q$ (the other implication is trivial). Since $\mathscr{A}_\sigma = \stcp\lbrace a_n \rbrace_{n \in \N} \subseteq \stcp{E} \subseteq Q$, $\mathscr{A}_\enlp = \mathscr{A}_\sigma \cup \mon \subseteq Q$. But then $\mathscr{A}_\enlp \subseteq \ga(Q)$ hence $\mon \subseteq \ga(Q)$
	\end{proof}
\end{prop}

\begin{cor}
	Given $\mon \in \M_{\aleph_0}(E)$, $Q \subseteq \enl{E}$ such that $\stcp{E} \subseteq Q$. If $\mon \subseteq Q$, then $\mon \cap \ga(Q) \neq \emptyset$.
	
	\begin{proof}
		Combine with corollary \ref{cor_countablemonadseq}.
	\end{proof}
\end{cor}

\subsection{I.o.G.-properties} \label{sec_iog}
Given a directed set $J$ and a family $(\gal_j)_{j \in J}$ in $\Ga(\enl{E})$, each having a base $\U_j$. Define:
\begin{displaymath}
	 \iog := \bigcap_{j \in \stcp{J}} \gal_j = \bigcap_{j \in \stcp{J}} \bigcup \stcp{\U_j};~~~ \gal := \ga(\iog) = \bigcup \stcp\big(\bigcap_{j \in J } \U_j \big);~~~ \I := \idl \iog = \idl \gal.
\end{displaymath}

We call $\iog$ an Intersection of Galaxies (I.o.G). For a standard set $V \subseteq E$, $\I$ is a set property that is characterized by
\begin{displaymath}
	 V \in \I \iff \enl{V} \subseteq \iog \iff \enl{V} \subseteq \gal
	\iff V \in \bigcap_{j \in J} \idl \gal_j.
\end{displaymath} 

We will consider three ways of extending this property to \textit{all} (including external) subsets of $E$. The $\enlp$-property is naturally defined by $\ee\I:= \cldownext\enl\I$. The $\cep$-extension (chromatic property) is $\ce\I:=\Pw(\iog)$.  Finally, the $\sep$-extension (standard property) is:
\begin{displaymath}
\se\I := \lbrace  A \subseteq \enl{E}: \forall j \in \stcp{J}: A \in \cldownext\enl\U_j \rbrace.
\end{displaymath} 
Each of them is (externally) downwardly closed.

\begin{prop} \label{prop_pmoniog}
Given $\mon \in \ps\M({E})$ with $\ps$base $\V$, we have
\begin{displaymath}
\mon \in \ce\I ~(\mon \subseteq \iog) \iff \mon \in \se\I \iff \forall j \in J: \exists V \in \V_\sigma, U \in \stcp\I: V \subseteq U.
\end{displaymath}

\begin{proof}
	The first equivalence follows from fixating $j \in \stcp{J}$ and applying the combined Cauchy principle on $\mon \subseteq \gal_j$. The second equivalence follows from fixating $j \in \stcp{J}$ and applying the Cauchy principle on $\mon \subseteq U$ for a certain $U \in \enl\U_j$.
\end{proof}
\end{prop}

\begin{cor}
	Given any set $A \subseteq \enl{E}$, then $A \in \ce\I ~(A \subseteq \iog)$ iff any internal set it contains is in $\se\I$. If $A$ is internal, then $A \in \ce\I ~(A \subseteq \iog) \iff A \in \se\I$. 
	
	Given $\gal' \in \ps\Ga(\enl{E})$ with $\ps$base $\U$, then $\gal' \in \ce\I ~(\gal' \subseteq \iog) \iff \U \subseteq \se\I$.
\end{cor}

Hence, we have
\begin{displaymath}
\iog = \bigcup_{(U_j) \in \prod_j \U_j} \bigcap_{j \in J} \enl{U_j};~~~~~ \gal = \bigcup_{(U_j) \in \prod_j \U_j} \enl\bigcap_{j \in J} {U_j}.
\end{displaymath}

\begin{prop}
Given $\gal' \in \Ga(\enl{E})$, we have 
\begin{displaymath}
 \gal' \in \ce\I ~(\gal' \subseteq \iog) \iff \gal' \in \ee\I  \iff \forall V \in \idl \gal': V \in \I \iff \gal' \subseteq \gal.
\end{displaymath}
\begin{proof}
We will first prove the second equivalence. Assume $\gal' \in \ee\I$ and take $V \in \idl \gal'$, then $\enl{V} \in \ee\I$, in fact in $\enl\I$ since it is internal which transfers to $V \in \I$. Assume otherwise the last statement is true. Then take any $N \in \enl(\idl \gal')_\infty$, because of the transfer principle $N \in \enl\I$, but then $\gal \in \ee\I$ because it is downward closed. The third equivalence is because of the Cauchy principle and the first and last statements are trivially equivalent.
\end{proof}
\end{prop}

\begin{prop} \label{prop_monioggal}
	Given $\mon \in \M(\enl{E})$, then the following are equivalent:	
	\begin{enumerate}[(i)]
		\item $\mon \subseteq \gal$;
		\item $\mon \in \ee\I$;
		\item $\fil \mon$ contains a set in $\I$;
		\item $\mon$ has a base in $\I$.
	\end{enumerate}
	\begin{proof}
		(i) $\iff$ (iii): This follows from the combined Cauchy principle.
		
		(ii) $\Rightarrow$ (iii): If $\mon \in \ee\I$, then each $U \in \enl(\fil \mon)_\infty$ is in $\ee\I$, in fact $\enl\I$ since they are internal; then because of the transfer principle there exists a $U \in \fil \mon \cap \I$. 
		
		(iii) $\Rightarrow$ (iv): If $\U$ is a base of $\mon$ and $V \in \I \cap \fil \mon$, then $\lbrace U \cap V : U \in \U \rbrace$ is a base of $\mon$ in $\I$.
		
		(iv) $\Rightarrow$ (ii): Trivial.
	\end{proof}
\end{prop}

In this case the distinction between $\enlp$-properties and $\sep$-properties are not merely non-standard tools, but have particular relevance in the standard world. For monads these carry over to two distinctive properties of filters, as will become clearer in concrete examples. 

\begin{lem}
	Given a ring $\Sigma$. If for all $j \in J: \gal_j$ is $\Sigma$-chromatic, then $\iog$ is $\Sigma$-chromatic.
	
\begin{proof}
	Since the $\Sigma$-chromatic sets are closed for intersection.
\end{proof}
\end{lem}

\subsection{Topology and convergence}

For any ring $\Sigma$, clearly if $U \in \Sigma$, then $\forall e \in U: \mo_\Sigma(e) \subseteq \enl{U}$. Suppose the opposite implication was also true, i.e. for any $U \subseteq E$:
\begin{displaymath}
	U \in \Sigma \iff \enl{U} = \bigcup_{e \in U} \mo_\Sigma(e) \textnormal{; i.e. } \forall e \in U: \exists V \in \fils{\Sigma} e: V \subseteq U;
\end{displaymath}
using the Cauchy principle. Note that this is strictly stronger than corollary \ref{cor_internalchromaticstandard} since the condition is generally weaker than being $\Sigma$-chromatic. In fact this statement is equivalent to $\Sigma$ being closed for the infinite union, i.e. a \textbf{topology}.

We will now work with a topological space $(E, \tau)$, i.e. $\tau$ is a ring that is closed for the infinite union, the \textbf{open sets}. Then $\overline{\tau}$ is closed for infinite intersection, the \textbf{closed sets}.

\begin{defn}
	For $e \in E$, the \textbf{topological monad} of $e$ is $\m_{[\tau]}(e) := \mo_\tau(e)$, $\tau$ is omitted when the topology is clear from the context. Then $\nbh_{[\tau]}(e) := \fil \m_{[\tau]}(e)$ is the filter of \textbf{neighborhoods} of $e$. Any base of $\m(e)$ is a base of neighborhoods of $e$.
\end{defn}
A set $U \subseteq E$ is thus open iff for any $e \in U$; $\m(e) \subseteq \enl{E}$, equivalently it contains a neighborhood of $e$. It is closed if it contains any $e$ such that $\m(e)\cap \enl{U} \neq \emptyset$. The \textbf{interior} of a set is its largest open subset and the \textbf{closure} of a set is its smallest closed superset. The $\enlp$-extensions of the open and closed properties are extended to external sets by letting $\ee\tau$ be the closure of $\enl\tau$ for infinite unions; equivalently the $\enl$closure of any $A \subseteq \enl{E}$ is
\begin{displaymath}
	\overline{A}^\tau := \lbrace e \in \enl{E}: \forall V \in \enl\nbh_\tau(e): A \cap V \neq \emptyset \rbrace.
\end{displaymath}
Since $\enl$interior and $\enl$closure are order preserving mappings, they have $\menlp$-extensions that coincide with the $\enlp$-extensions on monads and galaxies. Hence, a [$\ps$]monad or [$\ps$]galaxy is $\enl$closed iff it has [$\ps$]base of [$\enl$]closed sets and $\enl$open iff it has [$\ps$]base of [$\enl$]open sets. This implies that in fact for any $e \in \enl{E}$: $\overline{\mo(e)}^\tau = \mo_{\overline{\tau}}(e) =: \overline{\mo}^\tau(e)$. In the case of standard monads and galaxies, because of corollary \ref{cor_mongal_chrom} $\enl$open coincides with being $\tau$-chromatic and $\enl$closed coincides with being $\overline\tau$-chromatic.

\begin{defn}	
	A monad $\mon \in \M(\enl{E}) \setminus \lbrace \emptyset \rbrace$ \textbf{converges} towards a point $e \in E$, denoted $\mon \xrightarrow{\tau} e$ if $\mon \subseteq \m(e)$, i.e. each neighborhood of $e$ contains $\mon$.
	
	The \textbf{accumulation points} of $\mon$ are
	\begin{displaymath}
	\accpoints_\tau{\mon} := \ker \overline{\mon}^\tau.
	\end{displaymath}
\end{defn}

\begin{prop} \label{prop_accpointsconv}
	Given  $\mon \in \M(\enl{E})$:
	\begin{displaymath}
	\begin{split}
	\accpoints_\tau\mon  &= \lbrace e \in E: \m(e) \cap \mon \neq \emptyset \rbrace = \lbrace e \in E: \exists \mon' \in \M(\mon): \mon' \xrightarrow{\tau} e \rbrace.
	\end{split}
	\end{displaymath}
	
	\begin{proof}
		Assume $e \in \accpoints_\tau\mon$. Take any $U \in \enl\nbh(e)_\infty$; since $e \in \overline\mon$ there exists an $e' \in U \cap \mon$, hence $\mo(e') \subseteq \m(e) \cap \mon$. Now assume otherwise $e \in E$ and there exists an $e' 
		\in \m(e) \cap \mon$. Now recall that $\overline{\mon}$ has a base $\U$ of closed sets. Since for each $U \in \U$ we have $\m(e) \cap \enl{U} \ni e'$, the closure implies $e \in U$. Therefore $e \in \overline\mon$.
	\end{proof}
\end{prop}
\begin{cor} \label{cor_closedaccum}
	For $e \in \enl{E}$, $\accpoints_\tau\mo{(e)} = \lbrace e' \in E: \mo(e) \xrightarrow{\tau} e' \rbrace$.
\end{cor}
\begin{cor}
	For $V \subseteq E$, $\overline{V}^\tau = \accpoints_\tau{\enl{V}}$.
\end{cor}

\begin{lem} \label{lem_closedsep}
	Given $A \subseteq \enl{E}$. Then $\mo_\tau(A)$ is $\enl$closed iff for any closed $V \subseteq E$ such that $A \cap \enl{V} = \emptyset$, $A$ and $\enl{V}$ are $\tau$-separated.
	\begin{proof}
		Assume $\mo_\tau(A)$ is $\enl$closed and take a closed $V \subseteq E$. Then $\compl{\mo_\tau(A)}$ is an $\enl$open galaxy that contains $\enl{V}$, since it is $\tau$-chromatic it must also contain $\mo_\tau(\enl{V})$; hence $\mo_\tau(\enl{V}) \cap \mo_\tau(A) = \emptyset$. Assume otherwise that for any closed $V$, $A$ and $\enl{V}$ are $\tau$-separated and take an open $U$ such that $A \subseteq \enl{U}$. Then $\compl{\enl{U}}$ is closed and thus contained in an open set $W$ that is disjunctive from $\mo_\tau(A)$, i.e. $\mo_\tau(A) \subseteq \compl{W} \subseteq \enl{U}$. Hence $\mo_\tau(A)$ has a base of closed sets.
	\end{proof}
\end{lem}

A topology is \textbf{T1} if for any $e \in E$, $\lbrace e \rbrace$ is closed (i.e. $\ker \m(e) = \lbrace e \rbrace$); equivalently any standard point is $\tau$-distinguishable from any other standard point.
A topology is \textbf{Hausdorff (T2)} if each pair of standard points is $\tau$-separated; equivalently there exist no monads that converge to more than one point. A topology is \textbf{regular} if $\m(e)$ is $\enl$closed for any $e \in E$; equivalently (lemma \ref{lem_closedsep}) each singleton and closed set are $\tau$-separated. Any regular space is Hausdorff iff it is (T1), in that case the topology is \textbf{Hausdorff regular (T3)}.

If $V \subseteq E$, then the \textbf{induced topology} of $V$ is determined by the monads $\forall e \in V: \m_V(e) := \m_{E,\tau}(e) \cap \enl{V}$.

\begin{defn}
	The \textbf{near-standard points} of $\enl{E}$ are
	\begin{displaymath}
	\ns_\tau(\enl{E}) := \bigcup_{e \in E} \m_\tau(e) = \lbrace e \in \enl{E}: \accpoints_\tau{\mo(e)} \neq \emptyset \rbrace.
	\end{displaymath}	
	Since for any $e \in E$, $e' \in \m(e): \mo_\Sigma(e') \subseteq \m(e)$, $\ns(\enl{E})$ is $\tau$-chromatic.
	If the topology is Hausdorff then for each $e \in \ns_\tau(\enl{E})$, we define $\stt{\tau} e$ to be the unique element in $\accpoints_\tau{\mo(e)}$, i.e. $\mo(e) \xrightarrow{\tau} \st e$. 
\end{defn}

Hence, if the space is Hausdorff, corollary \ref{cor_closedaccum} implies that for $V \subseteq \enl{E}$:
\begin{displaymath}
	\overline{V}^\tau = \stt{\tau} \enl{V}.
\end{displaymath}

\begin{lem} \label{lem_stclosed}
	Suppose $E$ is Hausdorff. Given $\mon \in \ps\M(\enl{E})$. Then $\st \mon$ is closed.
	
	\begin{proof}
		Suppose $e \notin \st \mon$, then $\m(e) \subseteq \compl{ \mon}$. Because of the combined principle of Cauchy, there is an open neighborhood $U$ of $e$ such that $\mon \cap \enl{U} = \emptyset$. Hence $\st \mon \subseteq \st \enl(\compl{U}) = E \setminus U$ since its closed, hence $\m(e) \cap \st \mon = \emptyset$. 
	\end{proof}
\end{lem}

\begin{prop} \label{prop_nsiog}
	$\ns(\enl{E})$ is the I.o.G. of $\enl$open covering galaxies of $\enl{E}$.
	\begin{proof}
		It is clear that any $\enl$open covering galaxy must contain $\ns(\enl{E})$. Suppose then that there exists an element $e \in \enl{E}$ contained in every $\enl$open covering galaxy, then since $\compl{\overline{\mo(e)}}$ is open it cannot be covering, hence there exists an $e' \in \accpoints\mo(e)$.
	\end{proof}
\end{prop}

\begin{prop} \label{prop_relcompact}
	Given a set $V \subseteq E$. The following are equivalent:	
	\begin{enumerate}[(i)]
		\item $\enl{V} \subseteq \ns_\tau(\enl{E})$;
		\item $\forall \mon \in \M(\enl{V}), \mon \neq \emptyset: \accpoints_{\tau}\mon \neq \emptyset$.
		\item Each open cover of $E$ contains a finite subcover of $V$.
	\end{enumerate}
	
	\begin{proof}
		(i) $\iff$ (ii): From the definition of $\ns(\enl{E})$ and proposition \ref{prop_accpointsconv} follows that $\mon \cap \ns_\tau(\enl{E}) \neq \emptyset \iff \accpoints_\tau\mon \neq \emptyset$.
		
		(ii) $\iff$ (iii): Every open cover $\U$ is the subbase of an $\enl$open covering galaxy $\gal$. By proposition \ref{prop_nsiog}, $\enl{V} \subseteq \gal$, hence by the Cauchy principle there must be finite sets in $\U$ that cover $V$.
	\end{proof}
\end{prop}

\begin{defn}
	A set $V \subseteq E$ that satisfies the equivalent conditions of proposition \ref{prop_relcompact} is called \textbf{relatively compact}.\footnote{In non-regular spaces this definition is not equivalent with the definition of having compact closure (see further). Since relative compactness is the associated set property of an I.o.G., we can use section \ref{sec_iog} to define $\enl$relative-compact, $\se$relative-compact and $\ce$relative-compact on the basis of open covering galaxies and immediately get the implications for monads and galaxies. However, the approach with convergence and accumulation points is natural and sufficient.} $E$ is a compact space if $\enl{E} \subseteq \ns(\enl{E})$. A set $V \subseteq E$ is \textbf{compact} if $\enl{V} \subseteq \ns_V{(\enl{V})}$, i.e. the near-standard points for the space $V$ with the induced topology of $E$.
\end{defn}

A set that that is both closed and relatively compact, is always compact since in that case $\ns_V(\enl{V}) = \ns_E(\enl{V})$. If the topology is Hausdorff, a compact set set is always closed since for $e \in \ns_E(\enl{V})$: $\accpoints\mo(e) \cap V \neq \emptyset \iff \st e \in V$. This however does not imply that the closure of a relatively compact set is (relatively) compact.

\begin{lem} \label{lem_regclosed}
	If $E$ is regular, then any $\tau$-chromatic $A \subseteq \ns(\enl{E})$ is $\enl$closed.
	
	\begin{proof}
		Take $e \in \overline{A}$. Since $\mo_\tau(e)$ is $\enl$open, there exists an $e' \in \mo_\tau(e) \cap A$, or due to lemma \ref{lem_equivdis} $e \in \overline\mo(e')$. Since $A \subseteq \ns(\enl{E})$ there exists an $\hat{e} \in E$ such that $e' \in \m(\hat{e}) \subseteq A$. Since $E$ is regular $\m(\hat{e})$ is $\enl$closed, thus $e \in \overline\mo(e') \subseteq \m(\hat{e}) \subseteq A$.
	\end{proof}
\end{lem}

\begin{cor}
	If $E$ is regular, the closure of a relatively compact set is compact.
	\begin{proof}
		The lemma implies that $\ns(\enl{E})$ itself is $\enl$closed, hence the $\enl$closure of any subset of $\ns(\enl{E})$ must still be in $\ns(\enl{E})$.
	\end{proof}
\end{cor}

\begin{cor} \label{cor_reg_compact_closed}
	If $E$ is regular, any disjunctive pair of a relatively compact set and a closed set is $\tau$-separated.
	
	\begin{proof}
		Combine with lemma \ref{lem_closedsep}.
	\end{proof}
\end{cor}

\begin{defn}
	The \textbf{compact points} of $E$ are
	\begin{displaymath}
	\comp_\tau(\enl{E}) := \ga(\ns_\tau(\enl{E})),
	\end{displaymath} 	
	i.e. the galaxy of relatively compact sets. If $E$ is regular, it is also the galaxy of compact sets and $\comp_\tau(\enl{E})$ is $\enl$closed.
\end{defn}

\begin{prop} (Local property) \label{prop_lokaal}
	Given a galaxy $\gal$. Then the following statements are equivalent:
	\begin{enumerate}[(i)]
		\item $\ns(\enl{E}) \subseteq \gal$
		\item every convergent monad $\mon$ has a base $\U \subseteq \idl \gal$.
		\item every point has a neighborhood in $\U$;
	\end{enumerate}	
	The statements are always true for $\tau$-chromatic covering galaxies. If in general $\stcp{E} \subseteq \gal \subseteq \ns(\enl{E})$, they are implied by $\gal$ being $\tau$-chromatic.
	\begin{proof}
		(i) $\Rightarrow$ (ii): Given a base $\V$, because of the combined Cauchy principle it contains a set in $U \in \idl \gal$, then consider the base $\lbrace V \cap U: V \in \V \rbrace$;
		
		(ii) $\Rightarrow$ (iii): Trivial;
		
		(iii) $\Rightarrow$ (i): The assumption implies $\m(e) \subseteq \gal$ for any $e \in E$.
		
		For the final remark, $\tau$-chromatic is equivalent to $\enl$open for galaxies, so that this follows from proposition \ref{prop_nsiog}.
	\end{proof}
\end{prop}

\begin{exmp}
	$E$ is \textbf{locally compact} iff $\comp(\enl{E})$ is $\tau$-chromatic iff each near-standard point is compact iff each point has a compact neighborhood.
\end{exmp}

\begin{defn}
	Given a topological space $F$ and a function $\psi: \enl{E} \rightarrow \enl{F}$. Then $\psi$ is $\psi$ is \textbf{$\ce$continuous} if for any  $e \in E$, $\psi(\m(e)) \subseteq \m(\psi(e))$; $\psi$ is \textbf{$\se$continuous} if for any $e \in \stcp{E}$ and $V \in \stcp\nbh(\psi(e))$, there exists a $U \in \stcp\nbh(\psi(e))$ such that $\psi(U) \subseteq \psi(V)$. Finally, $\psi$ is \textbf{$\enl$continuous} if for any $e \in \enl{E}$ and $V \in \enl\nbh(\psi(e))$, there exists a $U \in \enl\nbh(\psi(e))$ such that $\psi(U) \subseteq \psi(V)$.
\end{defn}

The Cauchy principle implies that an \textit{internal} function $\psi$ is $\se$continuous iff it is $\ce$continuous. For standard functions all three properties coincide with continuity.

\section{Topological Vector Spaces}

\subsection{Linear monads}

We work with a $\K$-vector space $(E, +, \cdot)$, where $\K$ is either $\R$ or $\C$. We have the usual topology on $\K$, which means that:
\begin{displaymath}
\begin{gathered}
\forall \alpha \in \K: \m_\K(\alpha) = \lbrace \beta \in \enl{\K}: \forall n \in \N: |\alpha - \beta| < 1/n \rbrace;\\
\ns(\enl{\K}) = \enl\K_\f = \lbrace \alpha \in \enl{\K}: \exists n \in \stcp\N: |\alpha| \leq n \rbrace.
\end{gathered}
\end{displaymath}

Given $\mon \in \ps\M(\enl{E})$. Then for any $\mon' \in \ps\M(\enl{E})$, $\mon' + \mon$ is a $\ps$monad in $\enl$E and for any $\mon' \in \ps\M(\enl{\K})$, $\mon' \cdot \mon$ is a $\ps$monad in $E$ (due to corollary \ref{cor_imagemonad}).

The balanced and convex hulls of sets have non-standard equivalents applicable to external sets. Given $A \subseteq \enl{E}$, the \textbf{$\enl$balanced hull} is
\begin{displaymath}
\bal(A) := \bigcup_{\substack{\alpha \in \enl{\K}\\|\alpha| \leq 1}} \alpha A.
\end{displaymath}
 The \textbf{$\enl$convex hull} of a subset $A \subseteq \enl{E}$ is
\begin{displaymath}
\co(A) := \big\lbrace \sum_{i=1}^\omega a_i x_i : (x_i)_{i \leq \omega} \textnormal{ internal in }A,  (a_i)_{i \leq \omega} \textnormal{ internal in }\enl{[0,1]}, \sum_{i=1}^\omega a_i = 1 \big\rbrace.
\end{displaymath}	
 The $\enl$convex and $\enl$balanced hull is then
\begin{displaymath}
\cobal(A) := \co(\bal(A)) = \big\lbrace \sum_{i=1}^\omega a_i x_i : (x_i) \textnormal{ int. in }A,  (a_i) \textnormal{ int. in }\enl{\K}, \sum_{i=1}^\omega |a_i| \leq 1 \big\rbrace	
\end{displaymath}
and the $\enl$closed, $\enl$convex, $\enl$balanced hull $\overline{\cobal}(A) := \overline{\cobal(A)}$.

These maps coincide with the $\enl$-extension of their standard equivalent, so that for a standard $V$ we have $\enl({\bal}(V)) = \bal(\enl{V})$ and $\enl({\co}(V)) = \co(\enl{V})$. Furthermore, the two maps define set properties: $A$ is \textbf{$\enl$balanced} iff $A = \bal(A)$ and \textbf{$\enl$convex} iff $A = \co(A)$. Note that because of the transfer principle, an internal set $A$ is convex precisely iff for any two $e, e' \in A$ and $\alpha \in \enl[0,1]$: $\alpha e + (1 - \alpha) e \in A$, however this condition is not sufficient for external sets. Then, since the maps are order-preserving the maps and associated properties have $\menlp$-extensions that coincide with the $\enl$-extensions on $\ps$monads and $\ps$galaxies. Hence, a $[\ps]$monad is $\enl$balanced iff it has a $[\ps]$base existing out of $[\enl]$balanced sets and $\enl$convex iff it has a $[\ps]$base existing out of $[\enl]$convex sets.

In particular, sets that are closed for the $\enl\K_\f$-scalar multiplication are $\enl$balanced. In fact, a set is closed for $\enl\K_\f$-scalar multiplication iff it is closed for $\K$-scalar multiplication and is $\enl$balanced.

\begin{defn}
	Given two sets $A \subseteq \enl{E}$ and $B \subseteq \enl{E}$. Then
	\begin{displaymath}
	\begin{split}	
	&B \textbf{ $\se$absorbs } A \textnormal{ if }
\exists \alpha_0 \in \stcp\R_+: \forall \alpha \in \enl\K, |\alpha| \geq \alpha_0: A \subseteq \alpha B;\\
	&B \textbf{ $\enl$absorbs } A \textnormal{ if }
	\exists \alpha_0 \in \enl\R_+: \forall \alpha \in \enl\K, |\alpha| \geq \alpha_0: A \subseteq \alpha B;\\
	&B \textbf{ $\ce$absorbs } A \textnormal{ if } \m_\K(0) \cdot A \subseteq B.
	\end{split}
	\end{displaymath}
When $A$ is a singleton, we refer to it as a point. We call a set \textbf{$\se$absorbing} if it $\se$absorbs all points of $\stcp{E}$; \textbf{$\ce$absorbing} if it $\ce$absorbs $\stcp{E}$; \textbf{$\enl$absorbing} if it $\enl$absorbs all points of $\enl{E}$. 
\end{defn}

\begin{prop}
	Given $\mon \in \ps\M(\enl{E})$ and $\gal \in \ps\Ga(\enl{E})$. Then $\mon$ $\ce$absorbs $\gal$ iff any internal set containing $\mon$ $\se$absorbs any internal set contained in $\gal$.	
	\begin{proof}
		Given an internal set $A$, $\mon_A := \m_\K(0) \cdot A$ is a $\ps$monad with $\ps$base $(\lbrace \alpha A: |\alpha| \leq 1/n \rbrace)_{n \in \enl\N}$. For any internal set $B$, because of the Cauchy principle $\mon_A \subseteq B$ is equivalent with $\exists n \in N: \forall \alpha \in \K, |\alpha| > n:  \alpha A \subseteq B$, i.e. $B$ $\se$absorbs $A$. Then given $\U$ a $\ps$base for $\gal$ and $\V$ a $\ps$base for $\mon$:
		\begin{displaymath}
		\m_\K(0) \cdot \gal \subseteq \mon \iff \bigcup_{A \in \U} \m_A \subseteq \mon \iff \forall B \in \V: \bigcup_{A \in \U} \m_A \subseteq B.
		\end{displaymath} 
	\end{proof}
\end{prop}

\begin{cor}
	Given internal sets $A,B \subseteq \enl{E}$. Then $A$ $\se$absorbs $B$ iff $A$ $\ce$absorbs $B$.
\end{cor}

\begin{cor}
	Given $\mon \in \M(\enl{E})$ and $\gal \in \Ga(\enl{E})$. Then $\mon$ $\ce$absorbs $\gal$ iff $\mon$ $\enl$absorbs any set in $\idl\gal$.	
	\begin{proof}
		Given $\mon \in \M(\enl{E})$ and $V \in \idl\gal$. If $\m_\K(0) \cdot \enl{V} \subseteq \mon$, any set in $\fil \mon$ absorbs $V$. Transfer implies any set in $\enl(\fil \mon)_\infty$ $\enl$absorbs $\enl{V}$, hence $\mon$ $\enl$absorbs $\enl{V}$.  If otherwise $\mon$ $\enl$absorbs $\enl{V}$, any set in $\fil \mon$ $\enl$absorbs $\enl{V}$, by transfer it absorbs $V$ and therefore it $\ce$absorbs $\stcp{V}$.
	\end{proof}
\end{cor}

A topology $\tau$ is \textbf{compatible with the vector structure} if addition and multiplication are continuous for this topology. In this case $(E, \tau)$ is a \textbf{topological vector space (TVS)}.
	
\begin{prop} \label{prop_optelling_continu}
	Addition is continuous for $\tau$ iff
	\begin{enumerate}
		\item $\m_\tau(e) = \m_\tau(0) + e$ for each $e \in E$;
		\item $\m_\tau(0)$ is closed for addition.
	\end{enumerate}
	\begin{proof}
		Addition is continuous iff $\m(e_1) + \m(e_2) \subseteq \m(e_1 + e_2)$. Assume that is true, than (2) follows from the case $e_1 = -e_2 = 0$. If $e_1 = 0$ en $e_2 = e$ then $\m(0) + e \subseteq \m(e)$, but if $e_1 = e_2 = e$ then $\m(e) - e \subseteq \m(0)$, proving (1). Assuming otherwise (1) and (2) are true, then $\m(e_1) + \m(e_2) = e_1 + e_2 + \m(0) = \m(e_1 + e_2)$.
	\end{proof}
\end{prop}

\begin{defn} \label{def_mequivalentie}
	If addition on $E$ is continuous, we define for any $e \in \enl{E}$:
	\begin{displaymath}
	\m_\tau(e) := \m_\tau(0) + e.
	\end{displaymath}
	Note that for a non-standard $e$, $\m(e)$ is a $\ps$monad but not necessarily a standard monad. As these are equivalence classes, we denote for $e, e' \in \enl{E}$:
	\begin{displaymath}
	e \approx_{[\tau]} e' \iff e - e' \in \m(0).
	\end{displaymath}
\end{defn}

\begin{prop} \label{prop_vermenigvuldiging_continu}
	Suppose addition on $E_\tau$ is continuous. Then multiplication is continuous iff
	\begin{enumerate}
		\item $\m_\tau(0)$ is closed for $\enl\K_\f$-scalar multiplication, i.e. it is $\enl$balanced and closed for $\stcp{\K}$-scalar multiplication.
		\item $\m_\tau(0)$ is $\ce$absorbing, i.e. $\m_\K(0) \cdot \stcp{E} \subseteq \m_\tau(0)$.
	\end{enumerate}
	\begin{proof}
		Multiplication is continuous iff $\m_\K(\alpha) \cdot \m(e) \subseteq \m(\alpha e)$ for each $\alpha \in \K$ and $e \in E$. Assume that is true, then (1) follows from the case $e = 0$ and (2) from $\alpha = 0$. Assuming otherwise (1) and (2) are true, then $\m_\K(\alpha) \cdot \m(e) = (\alpha + \m_\K(0)) \cdot (e + \m(0)) = \alpha e + \m(0) = \m(\alpha e)$.
	\end{proof}
\end{prop}

\begin{cor}
	If $E_\tau$ is compatible with the vector structure, every neighborhood of $0$ is absorbing and contains a balanced neighborhood of $0$.
\end{cor}

\begin{defn}
	A monad $\mon \in \M(\enl{E})$ is \textbf{linear} if it is closed for $\enl\K_\f$-multiplication and it contains $\m_\K(0) \cdot \stcp{E}$.
	Any linear $\mon \in \M(\enl{E})$ \textit{defines} a topology $\tau$ on $E$ compatible with the topology by setting for $e \in E$:
	\begin{displaymath}
	\m_\tau(e) := \mon + e.
	\end{displaymath}	
\end{defn}

We will from now on assume that $E_\tau$ is a Hausdorff TVS, where $\Lambda_{[\tau]}$ is a base of $\m_\tau(0)$ existing out of absorbing, open and balanced sets. Given $\lambda \in \Lambda$ and $\alpha \in \K$, $\m(0) \subseteq \alpha \lambda$, then applying the Cauchy principle there exists a $\lambda' \in \Lambda$ such that $\lambda' \subseteq \alpha \lambda$. Furthermore, applying the Cauchy principle to $\m(0) + \m(0) \subseteq \lambda$ implies that for any $\lambda \in \Lambda$ there exists a $\lambda' \in \Lambda$ such that $\lambda' + \lambda' \subseteq \lambda$. 

Note that for $V \subseteq E$, $\overline{V} = \bigcap_{\lambda \in \Lambda} (V + \lambda)$, while for $A \subseteq \enl{E}$, $\overline{A} = \bigcap_{\lambda \in \enl\Lambda} (A + \lambda)$.

\begin{defn}
	A set $A \subseteq \enl{E}$ is \textbf{$\approx_\tau$-saturated} if $A = A + \m_\tau(0)$.
\end{defn}

Note that this is not equivalent to $\tau$-chromatic since for non-standard points $e' \in \enl{E}$, in general $\mo_\tau(e) \neq \m(e)$; they may be larger or smaller.

\begin{lem} \label{lem_sat_open}
	Any $\approx_\tau$-saturated set is $\enl$open and $\enl$closed.
	\begin{proof}
		Let $A$ be a $\approx_\tau$-saturated set and take $e \in A$. For any $\lambda \in \enl\Lambda_\infty$ we have $e + \lambda \subseteq A$, hence $A$ is $\enl$open.  Now take $e \in \overline{A}$. For any $\lambda \in \enl\Lambda_\infty$ we have $e' \in A$ such that $e - e' \in \lambda$, but then $e \in \m(e') \subseteq A$, hence $A$ is $\enl$closed.
	\end{proof}	
\end{lem}

\begin{cor}
	$\m_\tau(0)$ is $\enl$closed, i.e. $E$ is regular. 
\end{cor}

This means that $E_\tau$ has a closed base of neighborhoods of $0$, in particular $\lbrace \overline\lambda: \lambda \in \Lambda \rbrace$ is such a base. Furthermore for topological vector spaces T1=T2=T3, i.e. $E$ is Hausdorff iff $\ker \m_\tau(0) = \lbrace 0 \rbrace$.

\begin{lem} \label{lem_compplusgesloten}
	Given a closed $V \subseteq E$ and relatively compact $K \subseteq E$. Then $K + V$ is closed.	
	\begin{proof}
		Suppose $e \in \compl{(V + K)}$. Then $(e - V) \cap K = \emptyset$. Since $e - V$ is closed and $K$ relatively compact, because of regularity (corollary \ref{cor_reg_compact_closed}) $\mo_\tau(e - \enl{V}) \cap \mo_\tau(\enl{K}) = \emptyset$ and since $\m(e) - \enl{V} \subseteq \mo_\tau(e - \enl{V})$ this implies $\m(e) \cap \enl(K + V) = \emptyset$.
	\end{proof}
\end{lem}

Given $e \in \cobal(\stcp{E})$, a galaxy that has the convex, balanced hulls of finite sets as base. Then we have internal sequences $(e_i)_{i \leq \nu}$ in $\stcp{E}$ and $(\alpha_i)_{i \leq \nu}$ in $\enl\K$ such that $\sum_{i = 0}^\nu |\alpha_i| \leq 1$ and $e = \sum_{i = 0}^\nu \alpha_i e_i$. But if $\lbrace e_i\rbrace_{i \leq \nu }$ is an internal subset of $\stcp{E}$, due to the Cauchy principle it must be standard and finite, hence we may assume w.l.o.g that $\nu \in \stcp{\N}$. These observations lead to the following definition:

\begin{defn}	
	The \textbf{quasistandard points} are the $\enl{\K}_\f$-span of the standard points, i.e.
	\begin{displaymath}
	\qs(\enl{E}) := \cobal(\stcp{E}) = \lbrace \sum_{i = 0}^n \alpha_i e_i: n \in \N, (\alpha_i)_{i \leq n} \textnormal{ in } \enl{\K}_\f, (e_i)_{i \leq n} \textnormal{ in } \stcp{E}  \rbrace.
	\end{displaymath}
	Note that this galaxy is independent of the topology. 
\end{defn}

\subsection{Boundedness and precompactness}

\begin{defn}
	The set of \textbf{finite points} is the following I.o.G.  (cfr. section \ref{sec_iog}):
		\begin{displaymath}
	\fin_\tau(\enl{E}) := \bigcap_{\lambda \in \Lambda_\tau} \enl{\K}_\f \cdot \lambda = \bigcap_{\lambda \in \Lambda_\tau}  \bigcup_{\alpha_\lambda \in \R_+} \enl(\alpha_\lambda\lambda);
	\end{displaymath}
	i.e. $e \in \enl{E}$ is finite iff for any $\lambda \in \Lambda$ there exists an $\alpha \in \R_+$ such that $e \in \alpha\lambda$. The set is $\tau$-chromatic, $\tau$-saturated, closed for addition and $\enl{\K}_\f$-multiplication. The corresponding galaxy is the set of \textbf{bounded points}:
	\begin{displaymath}
	\bdd_\tau(\enl{E}) := \ga(\fin_\tau(\enl{E}))
	= \bigcup \stcp \lbrace \bigcap_{\lambda \in \lambda} \alpha_\lambda\lambda \rbrace_{(\alpha_\lambda)_\lambda \in \R_+^\Lambda} 
	\end{displaymath}	
	i.e. $e \in \enl{E}$ is bounded iff there exists a $(\alpha_\lambda)_{\lambda} \in \R_+^\Lambda$ such that $\forall \lambda \in \enl{\Lambda}: e \in \alpha_\lambda \lambda$.
\end{defn}

\begin{prop} \label{prop_charfin}
	A point $e \in \enl{E}$ is finite iff it is $\ce$absorbed by $\m_\tau(0)$.
	\begin{proof}
		Assume $e \in \fin_\tau(\enl{E})$ and $\epsilon \in \m_\K(0)$. For any $\lambda \in \Lambda$ there exists an $\alpha_\lambda \in \R_+$ such that $e \in \enl(\alpha_\lambda \lambda)$. Then $\epsilon e \in  \epsilon \alpha_\lambda(\enl\lambda) \subseteq \enl\lambda$ since $\lambda$ is balanced and $|\epsilon\alpha_\lambda| < 1$. Assume otherwise $e \notin \fin_\tau(\enl{E})$. There exists a $\lambda \in \Lambda$ such that $\forall \alpha \in \K$, $e \notin \enl(\alpha\lambda)$. Because of infinite overspill there exists an $\omega \in \enl\K_\infty$ such that $e \notin \omega(\enl\lambda)$, hence $\omega^{-1}e \notin \enl\lambda$.
	\end{proof}	
\end{prop}

The corresponding set property for $V \subseteq E$, $\enl{V} \subseteq \fin(\enl{E})$ is \textbf{bounded}. Following the notations of section \ref{sec_iog}; $A \subseteq \enl{E}$ is $\ce$bounded if $A \subseteq \fin(\enl{E})$ and [$\se|\enl$]bounded if for any $\lambda \in [\stcp|\enl]\Lambda$ there exists an $\alpha \in [\stcp|\enl]\R_+$ such that $A \subseteq \alpha \lambda$. As proven in the section, for internal sets and $\ps$monads, $\ce$boundedness is equivalent to $\se$boundedness, while for standard galaxies $\ce$boundedness is equivalent to $\enl$boundedness. Given a monad $\mon \in \M(\enl{E})$ we distinguish:
\begin{itemize}
	\item $\mon \subseteq \fin(\enl{E})$ is $\ce$bounded, $\se$bounded or \textbf{finite}, i.e. for any $\lambda \in \Lambda$ there exists an $\alpha \in \R_+$ and a $U$ in $\fil \mon$ such that $U \subseteq \alpha \lambda$; in light of proposition \ref{prop_charfin} a monad is finite iff when multiplied with a monad in $\enl\K$ that converges to $0$, the resulting monad converges to $0$;
	\item $\mon \subseteq \bdd(\enl{E})$ is $\enl$bounded or just \textbf{bounded}, i.e. $\fil \mon$ contains a bounded set and therefore a base existing out of bounded sets.
\end{itemize}
Hence, a set in $E$ is bounded iff each monad it contains is finite; a set in $\enl{E}$ is $\ce$bounded iff each monad it intersects is finite.

\begin{prop} \label{prop_cauchy}	
	Given $\mon \in \ps\M(\enl{E}) \setminus \lbrace \emptyset \rbrace$, the following statements are equivalent:	
	\begin{enumerate}[(i)]
		\item For each $e \in \mon$; $\mon \leq \m_\tau(e)$;
		\item there exists an $e \in \enl{E}$ such that $\mon \leq \m_\tau(e)$;
		\item for each $\lambda \in \Lambda_\tau$ there exists an internal (standard if $\mon$ is standard) $U_\lambda \supseteq \mon$ such that $U_\lambda - U_\lambda \subseteq \enl\lambda$;		
	\end{enumerate}
	
	\begin{proof}
		(i) $\Rightarrow$ (ii): Trivial.
		
		(ii) $\Rightarrow$ (iii): Given $\lambda \in \Lambda$, take $\lambda'$ such that $\lambda' + \lambda' \subseteq \lambda$. Because of the Cauchy principle there exists an internal $U_\lambda \supseteq \mon$ such that $U_\lambda \subseteq e + \enl\lambda'$. Then $U_\lambda - U_\lambda \subseteq \enl\lambda' + \enl\lambda' \subseteq \enl\lambda$.
		
		(iii) $\Rightarrow$ (i): Given $e \in \mon$ and take $\lambda \in \Lambda$. Because of the assumption there exists a $U_\lambda \supseteq \mon$ such that $U_\lambda - U_\lambda \subseteq \enl\lambda$, since $e \in U_\lambda$, $U_\lambda \subseteq e + \enl\lambda$.
	\end{proof}
\end{prop}

\begin{defn}
	A $[\ps]$monad that satisfies the equivalent conditions of \ref{prop_cauchy} is called \textbf{$[\ps]$Cauchy}.
\end{defn}

\begin{defn}
	The set of \textbf{prenearstandard points} is the following I.o.G.  (cfr. section \ref{sec_iog}):
	\begin{displaymath}
	\pns_\tau(\enl{E}) := \bigcap_{\lambda \in \Lambda_\tau}  \stcp{E} + \enl\lambda = \bigcap_{\lambda \in \Lambda_\tau} \bigcup_{e \in \stcp{E}} e + \enl\lambda;
	\end{displaymath}
	i.e. $e \in \enl{E}$ is prenearstandard iff for any $\lambda \in \Lambda$ there exists an $e' \in E$ such that $e - e' \in \lambda$. The set is $\tau$-chromatic, $\tau$-saturated, closed for addition and $\enl{\K}_\f$-multiplication. The corresponding galaxy is the set of \textbf{precompact points}:
	\begin{displaymath}
	\pcomp_\tau(\enl{E}):= \ga(\pns_\tau(\enl{E}))
	= \bigcup \stcp \lbrace \bigcap_{\lambda \in \lambda} S_\lambda + \lambda \rbrace_{(S_\lambda)_{\lambda} \in \mathcal{P}_\fin(E)^\Lambda}
	\end{displaymath}	
	i.e. $e \in \enl{E}$ is precompact iff there exists a $(S_\lambda)_{\lambda} \in \mathcal{P}_\fin(E)^\Lambda$ such that $\forall \lambda \in \enl{\Lambda}: e \in (\enl{S})_\lambda + \enl\lambda.$
\end{defn}

\begin{prop} \label{prop_pnsechtemonade}
	Given $e \in \enl{E}$. Then $e \in \pns_\tau(\enl{E})$ iff $\m_\tau(e) \in \M(\enl{E})$.	
	\begin{proof}
		Assume $e \in \pns(\enl{E})$. Then we can choose $(e_\lambda)_\lambda$ in $E$ such that for each $\lambda \in \Lambda$, $e - e_\lambda \in \lambda$. Then consider the monad $\mon$ with subbase $(e_{\lambda} + \lambda)_{\lambda \in \Lambda}$. Given $\lambda \in \Lambda$, take $\lambda'$ such that $\lambda' + \lambda' \subseteq \lambda$. For any $e' \in e_{\lambda'} + \lambda'$, $e - e' = (e - e_{\lambda'}) - (e' - e_{\lambda'}) \in \lambda' + \lambda' \subseteq \lambda$, so that $e_{\lambda'} + \lambda' \subseteq e + \lambda$; thus $\mon \subseteq \m(e)$. For any $e'' \in e + \lambda'$ similarly $e_\lambda - e'' \in \lambda$, so that $e + \lambda' \subseteq e_{\lambda} + \lambda$; thus also $\m(e) \supseteq \mon$. Hence $\m(e) = \mon \in \M(\enl{E})$.
		Assuming otherwise $\m(e) \in \M(\enl{E})$, then $\m(e)$ is cauchy. Given $\lambda$, there exists a $U \subseteq E$ such that $\m(e) \subseteq \enl{U}$ and $U - U \subseteq \lambda$. Picking any $e' \in U$, we have $e - e' \subseteq \lambda$.
	\end{proof}
\end{prop} 

The associated set property for $V \subseteq E, \enl{V} \subseteq \pns(\enl{E})$ is \textbf{precompact or totally bounded}. Following the notations of section \ref{sec_iog}, $A \subseteq \enl{E}$ is $\ce$precompact iff $A \subseteq \pns(\enl{E})$ and $[\se|\enl]$precompact iff for any $\lambda \in [\stcp|\enl]\Lambda$ there exists an $S \in [\stcp|\enl]\Pwfin(E)$ such that $A \subseteq S + \lambda$. As proven in the section, for internal sets and $\ps$monads $\ce$precompactness is equivalent to $\se$precompactness while for standard galaxies $\ce$precompactness is equivalent to $\enl$precompactness. Proposition \ref{prop_pnsechtemonade} implies that a standard monad is $\ce$precompact iff it is a union of Cauchy monads. Given a monad $\mon \in \M(\enl{E})$ we distinguish:
\begin{itemize}
	\item $\mon \subseteq \pns(\enl{E})$ is $\ce$precompact, $\se$precompact or \textbf{totally finite}, i.e. for any $\lambda \in \Lambda$ there exists a $U$ in $\fil \mon$ and a $S \in \Pwfin(E)$ such that $U \subseteq S + \lambda$. Each ultramonad it contains is Cauchy, i.e. for any $\lambda \in \Lambda$ there exists a $U$ in $\fil \mon$ such that $U - U \subseteq \lambda$;
	\item $\mon \subseteq \pcomp(\enl{E})$ are $\enl$precompact or just \textbf{precompact}, i.e. $\fil \mon$ contains a precompact set and therefore a base existing out of precompact sets.
\end{itemize}
Hence, a set in $E$ is precompact iff each ultramonad it contains is Cauchy, a set in $\enl{E}$ is $\ce$precompact iff each ultramonad it intersects is Cauchy.

\begin{rem}
	Note that because of proposition \ref{prop_sequentialmonads}, sequential monads are finite iff they are bounded and totally finite iff they are precompact. From its corollary follows that monads with countable cofinality that are finite, respectively totally finite always intersect with a bounded, respectively precompact set. This is the basis of many special properties of metrizables spaces where $\m_\tau(0)$ has countable cofinality.
\end{rem}

\begin{prop} \label{prop_lokaal2} (Local property for TVS) Given a covering $\gal \in \Ga(\fin_\tau(\enl{E}))$, closed for $\K$-scalar multiplication. The following statements are equivalent:
	\begin{enumerate}[(i)]
		\item $\m(0) \subseteq \gal$;
		\item $\fin_\tau(\enl{E}) = \gal$.	
		\item $\gal$ is $\tau$-chromatic;
		\item $\gal$ is $\approx_\tau$-saturated.
		\item There exists a neighborhood of $0$ in $\idl \gal$;
		\item Each finite monad has a base in $\idl \gal$.
	\end{enumerate}	
	\begin{proof}
		(i) $\Rightarrow$ (ii) : 
		Take $e \in \fin(\enl{E})$. Proposition \ref{prop_charfin} gives $\m_\K(0) \cdot e \subseteq \m(0) \subseteq \gal$, then there exists a $\alpha \in \R_+$ (combined principle of Cauchy) such that $\alpha e \in \gal$. Since $\gal$ is closed for $\K$-scalar multiplication, $e \in \gal$..
		
		(ii) $\Rightarrow$ (iii), (iv): Trivial.
		
		(iii), (iv) $\Rightarrow$ (i): Trivial.
		
		(i) $\iff$ (v): The combined Cauchy principle.
		
		(ii) $\iff$ (vi): Chromatic inclusion and the combined Cauchy principle.
	\end{proof}
\end{prop}

\begin{exmp} \label{ex_lokaalbegrensd}
	$E$ is \textbf{locally bounded}, i.e. $\m(0) \subseteq \bdd(\enl{E})$ iff $\fin(\enl{E}) = \bdd(\enl{E})$, in this case there exists a bounded $\lambda \in \Lambda$. For any $\lambda' \in \Lambda$ there exists an $\alpha \in \R_+$ such that $\alpha^{-1} \lambda \subseteq \lambda'$. Hence,
	\begin{displaymath}
	\m(0) = \bigcap_{\alpha \in \R_+}  \enl(\alpha^{-1}\lambda) = \m_\K(0) \cdot \enl\lambda.
	\end{displaymath}
	
	Looking back at the I.o.G. definition of $\fin(\enl{E})$ with galaxies of the form $\bigcup_{\alpha_\lambda \in \R_+} \enl(\alpha_\lambda\lambda)$, consider that for a bounded $\lambda$, this galaxy contains all the others so that $\fin(\enl{E})$ is indeed a single galaxy.
\end{exmp}

Since $\ns(\enl{E})$ is closed for (finite) addition and $\enl\K_\f$-scalar multiplication, $\cobal(\stcp{E}) \subseteq \ns(\enl{E})$.
Since $\stcp{E} \subseteq \pcomp(\enl{E})$ and the latter is $\approx_\tau$-saturated, $\ns(\enl{E}) \subseteq \pns(\enl{E})$. Given any $S \in \Pwfin(E)$ and $\lambda \in \Lambda$; there exists an $\alpha \in \R_+$ such that $S + \lambda \subseteq \alpha \lambda$, hence $\pns(\enl{E}) \subseteq \fin(\enl{E})$. Therefore we have
\begin{displaymath}
\begin{gathered}
\qs(\enl{E}) \subseteq \ns(\enl{E}) \subseteq \pns(\enl{E}) \subseteq \fin(\enl{E});\\
\qs(\enl{E}) \subseteq \comp(\enl{E}) \subseteq \pcomp(\enl{E}) \subseteq \bdd(\enl{E}).\\
\end{gathered}
\end{displaymath}
Hence, convex balanced hulls of finite sets are relatively compact, relatively compact sets are precompact and precompact sets are bounded. Convergent monads are Cauchy and Cauchy monads are finite. We may now derive special properties of (subspaces of) topological vector spaces from chromatic inclusions that do not generally hold, each with their equivalent characteristics:
\begin{itemize}  
	\item $\bdd(\enl{E}) \subseteq \pns(\enl{E})$. Bounded sets are precompact. Bounded monads are totally finite. Bounded ultramonads are Cauchy. \textbf{Pseudo-Henson-Moore}.	
	\item $\fin(\enl{E}) \subseteq \pns(\enl{E})$. Finite monads are totally finite. Finite ultramonads are Cauchy. \textbf{Henson-Moore}.\footnote{Named in \cite{stroyanlux} after \cite{hensonmoore1974}.} Implies Pseudo-Henson-Moore.  	
	\item $\pcomp(\enl{E}) \subseteq \ns(\enl{E})$. Precompact sets are compact. Precompact monads have an accumulation point. \textbf{Pseudo-complete}.
	\item $\pns(\enl{E}) \cap \bdd(\enl{E}) \subseteq \ns(\enl{E})$. Bounded closed sets are complete.\footnote{A set $V \subseteq E$ is complete if $\ns_V(\enl{V}) \subseteq \pns(\enl{V}).$ Closed subsets of a complete set are complete.} Bounded totally finite monads have an accumulation point. Bounded Cauchy monads are convergent. \textbf{Quasi-complete}. Implies pseudo-complete.
	\item $\pns(\enl{E}) \subseteq \ns(\enl{E})$. Totally finite monads have an accumulation point. Cauchy monads are convergent. \textbf{Complete}. Implies quasi-complete.
	\item $\bdd(\enl{E}) \subseteq \ns(\enl{E})$. Bounded sets are compact. Bounded monads have an accumulation point. Bounded ultramonads are convergent. \textbf{Heine-Borel}. Implies Quasi-complete and Pseudo-Henson-Moore. 
	\item $\fin(\enl{E}) \subseteq \ns(\enl{E})$. Finite monads have an accumulation point. Finite ultramonads converge.  \textbf{Henson-Moore-Complete}. Implies Heine-Borel, Henson-Moore and complete.
\end{itemize}

We immediately see that a set is compact iff it is precompact and complete. 

\begin{prop} \label{prop_lokaalprecompactiseindigdim}
		The following statements are equivalent:	
		\begin{enumerate}[(i)]
			\item $\m_\tau(0) \subseteq \comp_\tau(\enl{E})$, i.e. $E$ is locally compact;
			\item $\m_\tau(0) \subseteq \pcomp_\tau(\enl{E})$, i.e. $E$ is locally precompact.
			\item $\m_\tau(0) \subseteq \qs(\enl{E})$;
			\item $E$ is finite dimensional;			
		\end{enumerate}	
		In this case $\fin(\enl{E}) = \bdd(\enl{E}) = \pns(\enl{E}) = \pcomp(\enl{E}) = \ns(\enl{E}) = \comp(\enl{E}) = \qs(\enl{E})$. In particular, $E$ is  Henson-Moore-Complete.
		
		\begin{proof}
			(i) $\Rightarrow$ (ii): Trivial.
			
			(ii) $\Rightarrow$ (iii): Let $\lambda$ be a precompact neighborhood of $0$. Hence, there exists an $S \in \Pwfin(E)$ so that recursively
			\begin{displaymath}
			\begin{split}
			\lambda &\subseteq S + \frac{\lambda}{2} \subseteq S + \frac{1}{2}(S + \frac{\lambda}{2}) \subseteq \cobal(S) + \frac{\lambda}{4} 
			\\&\subseteq .... \subseteq \cobal(S) + \frac{\lambda}{2^n} ~~(\forall n \in \N) \Rightarrow \lambda \subseteq \overline\cobal(S).
			\end{split}
			\end{displaymath}
						
			(iii) $\iff$ (iv): The galaxy with base $\lbrace \spann(S) \rbrace_{S \in \Pwfin(E)}$ (finite dimensional subspaces) is equal to $\enl\K \cdot \qs(\enl{E})$, which is equal to $\enl{E}$ if $E$ is finite-dimensional. If $E$ is not finite-dimensional, then for any finite $S \subseteq E$ there exists an $e \in \compl{\spann(S)}$, then for $\epsilon \in \m_\K(0)$, we have $\epsilon e \in \m(0)$ thus $\m(0) \nsubseteq \enl\spann(S)$. Since $S$ was arbitrary, $\m(0) \nsubseteq \enl\K \cdot \qs(\enl{E})$.
			
			(iii) $\Rightarrow$ (i): Since $\qs(\enl{E})$ is covering and closed for addition, the assumption implies $\ns(\enl{E}) = \qs(\enl{E})$. Since the latter is a galaxy it must be equal to $\comp(\enl{E})$.				
		\end{proof}
	\end{prop}

\begin{cor} 
	The only topology on Hausdorff finite dimensional spaces is the Euclidean topology.	
	\begin{proof}	
		Indeed, if $E$ is finite-dimensional, $\m(0) = \m(0) \cap \qs(\enl{E})$, the $\m_\K(0)$-span of $\stcp{E}$, which is defined independently of the topology.
	\end{proof}
\end{cor}

\begin{cor} 
	$\qs(\enl{E}) = \overline{\qs(\enl{E})}^\tau$, regardless of $\tau$ (convex balanced hulls of finite sets are closed, hence compact).
\end{cor}

\subsection{Completion}

We can identify $E$ with the quotient space $\ns(\enl{E})/\m(0)$, then $\st \cdot$ is identified with the projection map. Consider now $\hat{E} := \pns(\enl{E})/\m(0)$ and let  $\stcom \cdot$ be the projection map for this quotient space. Then $E$ is identified with a subspace of $\hat{E}$ and $\st \cdot$ is the restriction of $\stcom \cdot$ to $\ns(\enl{E})$. Now, consider the map 
\begin{displaymath}
	c: \Pw(E) \rightarrow \Pw(\hat{E}): V \rightarrow \stcom \enl{V}.
\end{displaymath}
Since $c$ is order-preserving, $\enl{c}$ is $\menlp$-extendable; so let
\begin{displaymath}
	\m_{\hat{E}}(0) := \menl{c}(\m_E(0)).
\end{displaymath}
Since projections are linear, $c$ is linear. Then it is clear that $\m_{\hat{E}}(0)$ must be closed for addition and $\enl\K_\f$-scalar multiplication. Furthermore, since $c(E) = \hat{E}$, we find for any $\epsilon \in \m_\K(0)$ that
\begin{displaymath}
\epsilon \cdot \stcp{\hat{E}} = \epsilon \cdot \menl{c}(\stcp E) = \menl{c}(\epsilon \cdot \stcp{E}) \subseteq \menl{c}(\m_E(0)) = \m_{\hat{E}}(0).
\end{displaymath}
Hence, the monad is linear and defines a topology $\hat{\tau}$ compatible with the vector structure. Since any $e \in \pns(\enl{E})$ is $\tau$-separated from the origin, there exists a $\lambda \in \Lambda$ that is $\tau$-separated from $e$ which implies that $\stcom e$ is not in $\m_{\hat{E}}(0)$, hence $\hat{\tau}$ is Hausdorff. Now for any $V \subseteq E$, $c(V) \cap E = \st \enl{V} = \overline{V}$, hence  $\m_{\hat{E}}(0) \cap E = \m_E(0)$, i.e. $\hat{\tau}$ induces $\tau$ on $E$. 
Suppose $e \in \pns(\enl\hat{E})$, then for any $\lambda \in \Lambda$, $e \in \stcp\hat{E} + c(\lambda) = \menl{c}(\stcp{E} + \lambda)$. Applying the Cauchy principle for each $\lambda$, there exist $(e_\lambda)_{e \in \Lambda}$ in $E$ such that $e \in  \bigcap_{\lambda} c(e_\lambda + \lambda) = \menl{c}(\m_E(e'))= \m_{\hat{E}}(e')$ for some $e' \in \pns(\enl{E})$. Then $e \approx \stcom e'$ so that $e \in \ns(\enl\hat{E})$. In other words, $\hat{E}$ is complete. Then $\stcomt{\tau} \cdot$ coincides with $\stt{\hat\tau} \cdot$ on $\enl{E}$ and $c$ coincides on $E$ with closure in $\hat{E}$, but then $\overline{E}^{\hat\tau} = \hat{E}$. Finally, it is clear that $\hat{E}$ is embeddable in any complete space containing $E$, so that $\hat{E}$ is (up to an isomorphism) the unique \textbf{completion} of $E$.

We can use the $\stcom \cdot$ map to generate weaker forms of completions: $\stcom \bdd(\enl{E}) \cap \pns(\enl{E})$ is the quasi-completion; $\stcom \pcomp(\enl{E})$ the pseudo-completion.

\begin{lem}
	Given an internal $A \subseteq E$ such that $A$ is $\se$finite. Then $\st A$ is finite.	
	\begin{proof}
		There exist $(\gamma_\lambda)_\lambda \in \K^\Lambda$ such that $A \subseteq \bigcap_{\lambda \in \stcp\Lambda} \enl(\gamma_\lambda \lambda)$. Then $\st A \\ \subseteq \st \bigcap_{\lambda \in \stcp\Lambda} \enl(\gamma_\lambda \lambda) \subseteq \bigcap_{\lambda \in \stcp\Lambda} \gamma_\lambda (\st \enl\lambda) = \bigcap_{\lambda \in \Lambda} \gamma_\lambda \overline\lambda$.
	\end{proof}
\end{lem}

\begin{lem} \label{lem_sttotaalbegrensd}
	Given an internal $A \subseteq E$ such that $A$ is $\se$precompact. Then $\st A$ is precompact.	
	\begin{proof}
		There exist $(S_\lambda)_\lambda \in \mathcal{P}_\fin(E)^\Lambda$ such that $A \subseteq \bigcap_{\lambda \in \stcp\Lambda} S_\lambda + \enl\lambda$. Then $\st A \subseteq  \bigcap_{\lambda \in \stcp\Lambda} S_\lambda + \st\enl\lambda = \bigcap_{\lambda \in \Lambda} S_\lambda + \overline\lambda$.
	\end{proof}
\end{lem}

The two lemma's above may be applied to $\stcom \cdot = \stt{\hat{\tau}} \cdot$ as well. In the case that $A$ is $\se$precompact, $\stcom A$ is compact.

\begin{lem} \label{lem_st_comp}
	Given an internal $A \subseteq \ns(\enl{E})$. Then $\st A = \stcom A$ is compact.	
	\begin{proof}
		Since $A$ is $\se$precompact, lemma \ref{lem_sttotaalbegrensd} implies that $\st A$ is precompact. Furthermore, $\st A$ is closed by lemma \ref{lem_stclosed}.  Since $\st \cdot = \stcom \cdot$ on $\ns(\enl{E})$, $\st A$ is in fact closed in $\hat{E}$ and therefore complete. 
	\end{proof}
\end{lem}

\subsection{Locally convex spaces} \label{sec_coco}

$E$ is \textbf{locally convex} iff $\m_\tau(0)$ is $\enl$convex. In this case we may assume that $\Lambda$ is a base existing out of absorbing, open, balanced and convex sets. Furthermore, we will in that case interpret $\Lambda$ as a family of continuous semi-norms, using the Minkowski-function
\begin{displaymath}
\|e\|_\lambda := \inf \lbrace r \in \R_+: e \in r \lambda \rbrace.
\end{displaymath}
such that the open (respectively closed) unit ball of the seminorm is precisely $\lambda$ (respectively $\overline\lambda$). Then $e \approx_\tau 0 \iff \forall \lambda \in \Lambda_\tau: \|e\|_\lambda \approx 0$; a single norm is sufficient (i.e. the space is \textbf{normable}, i.e. induced by a norm) iff it is locally bounded.

We will assume that $E$ is locally convex in this section.

\begin{prop} \label{prop_cobalfin}
	 $\overline\cobal(\fin_\tau(\enl{E})) = \fin_\tau(\enl{E}).$	
	\begin{proof}
		Since $\overline\cobal(\enl\K_\f \cdot \overline\lambda) = \enl\K_\f \cdot \overline\lambda$, this follows from lemma \ref{lem_hullintersect}.
	\end{proof}	
\end{prop}

\begin{cor}
	The closed convex balanced hull of a bounded set is bounded, hence
	$\overline\cobal({\bdd_\tau(\enl{E})}) = \bdd_\tau(\enl{E})$.	
\end{cor}

\begin{prop} \label{prop_cobalpns}
	$\overline\cobal(\pns_\tau(\enl{E})) = \pns_\tau(\enl{E}).$
	\begin{proof}
		Since $\stcp{E} \subseteq \qs(\enl{E}) \subseteq \pns(\enl{E})$, 
		\begin{displaymath}
			\pns(\enl{E}) = \bigcap_{\lambda \in \Lambda} \qs(\enl{E}) + \overline\lambda.
		\end{displaymath}
		Now $\overline\cobal(\qs(\enl{E}) + \overline\lambda) = \qs(\enl{E}) + \overline\lambda$. Indeed, the sum of two $\enl$convex, $\enl$balanced sets is $\enl$convex and $\enl$balanced and because of lemma \ref{lem_compplusgesloten} this galaxy is also $\enl$closed. Hence, the result follows from lemma \ref{lem_hullintersect}.
	\end{proof}
\end{prop}

\begin{cor}
The closed convex balanced hull of a precompact set is precompact, hence
$\overline\cobal({\pcomp_\tau(\enl{E})}) = \pcomp_\tau(\enl{E})$.	
\end{cor}

\begin{prop} \label{prop_cobalnscomp}	
	$\overline\cobal(\ns_\tau(\enl{E})) = \cobal(\comp_\tau(\enl{E})) + \m(0)$.
	\begin{proof}
		Take $e \in \enl{\cobal}(\ns(\enl{E}))$. There exist internal sequences $(e_n)_{n \leq \omega}$ in $\ns(\enl{E})$ en $(\alpha_n)_{n \leq \omega}$ in $\enl\K$ such that $\sum_{i = 0}^\omega |\alpha_i| \leq 1$ and $e = \sum_{i = 0}^\omega \alpha_i e_i$. Let $U := \st \lbrace e_n \rbrace_{n \leq \omega}$, a compact set due to lemma \ref{lem_st_comp}. Now $\lbrace e_n \rbrace_{n \leq \omega}$ is $\se$precompact, so that because of infinite overspill there exists a $\lambda \in \enl\Lambda_\infty$ and $S_\lambda \in \enl\mathcal{P}_\fin(V)$ such that $\lbrace e_n \rbrace_{n \leq \omega} \subseteq S_\lambda + \lambda$. For each $e_i$ we find an $e'_i \in S_\lambda \subseteq \enl{V}$ such that $e_i - e'_i \in \lambda$, hence $e_i \approx e'_i$. Let $e' := \sum_{i = 0}^\omega \alpha_i (e'_i)$. Then $e' \in \enl\cobal(V) \subseteq \cobal(\comp(\enl{E}))$.For arbitrary $\lambda \in \Lambda$ we find
		\begin{displaymath}
		\|e - e'\|_\lambda \leq \sum_{i = 0}^\omega |\alpha_i| \|e_i - e'_i\|_\lambda
		\leq (\sum_{i = 0}^\omega |\alpha_i|) \max \|e_i - e'_i\|_\lambda \approx 0
		\end{displaymath}
		so that $e' \approx e$, hence $e \in \enl\cobal(V) + \m(0) \subseteq \cobal(\comp_\tau(\enl{E})) + \m(0)$.
		
		Now since $\m(0) \subseteq \cobal(\ns_\tau(\enl{E}))$ and the latter is closed for addition, $\cobal(\ns_\tau(\enl{E}))$ is $\approx_\tau$-saturated, and thus also $\enl$closed; from which the result follows.
	\end{proof}
\end{prop}

\begin{cor} \label{cor_cobalns}	
The closed convex balanced hulls of compact sets are compact ($E$ has \textbf{the convex convex compactness property}) iff $\ns(\enl{E})$ is $\enl$convex iff $\comp(\enl{E})$ is $\enl$convex.
\begin{proof}
	Note that $\ns(\enl{E})$ and $\comp(\enl{E})$ are always $\enl$balanced and $\enl$closed. Then if $\comp(\enl{E}) = \overline\cobal(\comp(\enl{E}))$; $\overline\cobal(\ns(\enl{E})) = \cobal(\comp(\enl{E}) + \m(0)) = \overline\cobal(\comp(\enl{E})) + \m(0) = \ns(\enl{E})$. On the other hand, if $\ns(\enl{E}) = \cobal(\ns(\enl{E}))$, convex hulls of compact sets are clearly compact.
\end{proof}
\end{cor}

The convex compactness property may be understood as a weak form of completeness. The smallest subspace of $\hat{E}$ containing $E$ that has the convex compactness property is $\hat{E}_c = \stcom \overline\cobal(\ns(\enl{E}))$. Then for a set $V \subseteq E$, we have $\enl{V} \subseteq \cobal(\ns(\enl{E}))$ iff all of its monads converge in $\hat{E}_c$. Furthermore, since $\overline\cobal(\comp(\enl{E})) \subseteq \pcomp(\enl{E})$, pseudo-completeness (and quasi-completeness) are sufficient conditions for the convex compactness property.

In general it is not so that sets contain a greatest balanced, convex subset (for the  $\subseteq$ order). In the case of a galaxy that is closed for addition there appears to be such a thing as a greatest balanced, convex subgalaxy.

\begin{lem}
	Given $\gal \in \Ga({E})$ that is closed for addition. Given two balanced convex sets $V, W \in \idl \gal$; $\cobal(\enl{V} \cup \enl{W}) \subseteq \gal$.	
	\begin{proof}
		Since $V$ and $W$ are balanced and convex, $\cobal(V \cup W) \subseteq V + W$. Since $\gal$ is closed for addition we must have $\enl\cobal(V \cup W) \subseteq \enl{V} + \enl{W} \subseteq \gal$.
	\end{proof}
\end{lem}

Hence, there exists a galaxy $\gal' \subseteq \gal$ that takes as its base all  $\enl$convex, $\enl$balanced sets of $\idl \gal$. Then $\gal'$ must contain any $\enl$convex, $\enl$balanced subgalaxy of $\gal$. If $\gal$ is $\enl$closed, this works for $\enl$convex, $\enl$balanced, $\enl$closed sets. In particular for compact points, the following definition makes sense:
\begin{defn}
	The \textbf{convex-compact points} is the galaxy of points in $\enl{E}$ which are contained in a convex, balanced, compact subset of $E$. We denote:
	\begin{displaymath}
	\coco(\enl{E}) := \bigcup \lbrace \enl{V}: \enl{V} \subseteq \ns(\enl{E}), \overline\cobal({V}) = V \rbrace.
	\end{displaymath}
\end{defn}

Note that $\coco(\enl{E}) = \comp(\enl{E})$ iff $E$ has the convex compactness property.

\subsection{Linear maps}

In this section we will work with a second topological vector space $(G, \tau')$ where $\m_{G, \tau'}(0)$ has a base $\Kappa$ with the same assumptions as $\Lambda$.

For linear maps, continuity is determined in the origin. This means that a linear map $\psi: \enl{E} \rightarrow \enl{G}$ is $\ce$continuous iff $\psi(\m_E(0)) \subseteq \m_G(0)$, as can easily be proven using linearity. There appears to be a second characterization:
\begin{prop} \label{prop_contbweind}
	Given a linear map $\psi: \enl{E} \rightarrow \enl{G}$. Then $\psi$ is $\ce$continuous iff $\psi(\fin(\enl{E})) \subseteq \fin(\enl{G})$.
	
	\begin{proof}	
		First assume $\psi$ is $\ce$continuous. Then since
		 \begin{displaymath}
		\m_\K(0) \cdot \psi(\fin(\enl{E})) = \psi(\m_\K(0) \cdot \fin(\enl{E})) \subseteq \psi(\m_E(0)) \subseteq \m_G(0)
		\end{displaymath} 
		this follows from proposition \ref{prop_charfin}. Now assume otherwise that $\psi(\fin(\enl{E})) \subseteq \fin(\enl{G})$. Since $\m_\K(0) \cdot \m_E(0)$ is the monad of neighborhoods of $0$ in $E$ multiplied with neighborhoods of $0$ in $\K$, it is in fact equal to $\m_E(0)$. Then, using proposition \ref{prop_charfin} again,
		 \begin{displaymath}
		\begin{split}
		\psi(\m_E(0)) &= \psi(\m_\K(0) \cdot \m_E(0)) = \m_\K(0) \cdot \psi(\m_E(0)) \\&\subseteq \m_\K(0) \cdot \psi(\fin(\enl{E})) \subseteq \m_\K(0) \cdot \fin(\enl{G}) \subseteq \m_G(0)
		\end{split}
		\end{displaymath} 	
	\end{proof}
\end{prop}

This means that $\ce$continuity (including $\se$continuity for internal maps and continuity for standard maps) is equivalent with mapping finite monads onto finite monads. For internal maps, this implies that bounded sets are mapped onto bounded sets, i.e. $\psi(\bdd(\enl{E})) \subseteq \bdd(\enl{G})$ (since galaxies are mapped onto galaxies). The latter however does not in general imply $\ce$continuity (it does so obviously when  $\bdd(\enl{G}) = \fin(\enl{G})$, i.e. the space is locally bounded).

\begin{prop} \label{prop_linafbmonad}
	Given a linear map $\theta: E \rightarrow G$. Then $\theta(\m(0))$ is a linear monad for the subspace $\theta(E)$ of $G$.
	
	\begin{proof}
		 Because of corollary \ref{cor_imagemonad}, $\theta(\m(0))$ is a monad in $G$ . Due to linearity of $\theta$, $\theta(\m(0))$ is closed for addition and $\K_\f$-scalar multiplication. We have $\ce$absorption since $\m_\K(0) \cdot \theta(E) = \theta(\m(0) \cdot E) \subseteq \theta(\m(0))$.
	\end{proof}
\end{prop}

Hence, a surjective linear map $E \rightarrow G$ defines a topology on $G$. If $G$ is a TVS, $\theta$ is a linear isomorphism and $\theta(\m_E(0)) = \m_F(0)$, i.e. $\theta$ is a homeomorfism, then $\theta$ is a \textbf{TVS-isomorfism}. 

Given a subspace $H \subseteq E$, and $\pi_H: E \rightarrow E/H$ is the projection to quotient space $E/H$, then $\pi_H(\m_E(0))$ is a linear monad of $E/H$.	Then $\m_{E/H}(0) := \pi_H(\m_E(0))$ defines a natural topology on $E/H$. Under this topology $\pi_H$ is open and continuous. If $H$ is closed, then $\pi_H(\compl{H}) = E/H \setminus \lbrace 0 \rbrace$ is open and $E/H$ is Hausdorff.	

\begin{defn}
	A family $\lbrace \psi_j \rbrace_{j \in J}$ of linear maps $E \rightarrow G$ is \textbf{equicontinuous} if for each $\kappa \in \Kappa$ there exists a $\lambda \in \Lambda$ such that $\psi_j(\lambda) \subseteq \kappa$ for all $j \in J$.
\end{defn}

\begin{prop} \label{prop_equicontinu}
	A family $\lbrace \psi_j \rbrace_{j \in J}$ of linear maps $E \rightarrow G$ is \textbf{equicontinuous} iff for each $j \in \enl{J}$, $\psi_j$ is $\ce$continuous (i.e. $\se$continuous).		
	\begin{proof} 
		Suppose that for each $j \in \enl{J}$, $\psi_j$ is $\ce$continuous. Take any $\kappa \in \Kappa$. Then 
		\begin{displaymath}
		\m_E(0) \subseteq \bigcap_{j \in \enl{J}} \psi_j^{-1}(\enl{\kappa}).
		\end{displaymath}
		so that due to the Cauchy principle there exists a $\lambda \in \Lambda$ such that $\enl{\lambda} \subseteq \psi_j^{-1}(\enl{\kappa})$ for each $j \in J$. Assume otherwise that the family is equicontinuous and take $j \in \enl{J}$. There exist for each $\kappa \in \Kappa$ a $\lambda_\kappa \in \Lambda$ such that due to the transferring principle $\psi_j(\lambda_\kappa) \subseteq \kappa$, hence
		\begin{displaymath}
		\psi_j(\m_E(0)) \subseteq \psi_j(\bigcap_{\kappa \in \Kappa}\enl{\lambda_\kappa}) \subseteq \bigcap_{\kappa \in \Kappa} \psi_j(\enl{\lambda_\kappa}) \subseteq \bigcap_{\kappa \in \Kappa} \enl\kappa = \m_G(0).
		\end{displaymath}
	\end{proof}
\end{prop}

\section{Duality Theory}

\subsection{Dual pairs and poles}

The following lemma will be assumed without proof; it is Hahn-Banach for finite dimensions and can be proven inductively by a well known argument:
\begin{lem} \label{lem_hb}
	Given an finite dimensional space $E$ with euclidean norm $\|.\|$, a subspace $H$ and a linear functional $\psi$ on $H$ with $\|\psi\| = \sup_{\|e\| \leq 1} \psi(e) \leq 1$. There exists a linear functional $\hat{\psi}$ on $E$ such that $\restr{\hat{\psi}}{H} = \psi$ and $\|\hat{\psi}\| = 1$.
\end{lem}

\begin{thm} (Hahn Banach)
	Given Hausdorff TVS $E$, a subspace $H$, a closed absorbing balanced convex set $D$ and a linear functional $\psi$ on $H$ such that $\sup_{e \in D \cap H} |\psi(e)| \leq 1$. There exists a linear functional $\hat{\psi}$ on $E$ such that $\restr{\hat{\psi}}{H} = \psi$ and $\sup_{e \in D} |\psi(e)| \leq 1$.
	\begin{proof}
		Consider the family $\mathcal{R}$ of finite dimensional subspaces of $E$; directed for $\subseteq$. For any $R \in \mathcal{R}$, since $E$ is Hausdorff the induced topology is Euclidean. Now take any $R \in \enl\mathcal{R}_\infty$, then $\stcp{E} \subseteq R$. Let $D_R := \enl{D} \cap R$, a $\enl$balanced, $\enl$convex and $\enl$closed set. Because of the transfer principle,  $D_R$ is of the form $\lbrace e \in R: \|e\| \leq 1 \rbrace$ for some $\enl$Euclidean\footnote{By this we mean $\|\sum_{j = 0}^\omega \alpha_j e_j\| = \sqrt{\sum_{j \leq \omega} \alpha_j^2}$ for a hyperfinite base $\lbrace e_j \rbrace_{j \leq \omega}$ of $R$.} norm. Apply lemma \ref{lem_hb}, using the transfer principle, on $R$ with subspace $H_R := \enl{H} \cap R$ and $\psi_R := \restr{\enl\psi}{H_R}$. Then there is a linear functional $\hat{\psi_R}$ on $R$ such that $\sup_{e \in D_R} |\hat{\psi_R}(e)| \leq 1$. Since $D$ is absorbing, for any $e \in \stcp{E}$ there exists an $\alpha \in \R_+$ such that $e \in \stcp{(\alpha D)} \subseteq \alpha D_R$, hence $\hat{\psi_R}(e) \leq \alpha \sup_{e \in D}  \hat{\psi_R}(D_R) \leq \alpha$. Hence $\hat{\psi_R}(e)$ is finite and we can set $\hat{\psi}(e) := \st \hat{\psi_R}(e).$
	\end{proof}
\end{thm}

\begin{cor} \label{cor_hbsep}
	Given a Hausdorff TVS $E$, a closed absorbing balanced convex $D \subseteq E$ and a point $x \in \compl{D}$. There exists a linear functional $\psi$ on $E$ such that $\sup_{e \in D} |\psi(e)| \leq 1$ and $|\psi(x)| > 1$.
	\begin{proof}
		Let $H = \K \cdot x$. Then there exists an $r \in \left] 0, 1\right[$ such that $H \cap D = \lbrace \alpha x: |\alpha| \leq r \rbrace $. Let $\psi$ be the extension on $E$ of the linear functional on $H$ defined by $\psi(\alpha x) := \alpha r^{-1}$.
	\end{proof}
\end{cor}

	A \textbf{dual pair} $(E, F)$ are two vector spaces $E$ and $F$ over a field $\mathbb{F}$ and a bilinear mapping $\langle , \rangle : E \times F \rightarrow \mathbb{F}$ such that both spaces are separated by each other:
	\begin{displaymath}
	\forall e \in E \setminus \lbrace 0 \rbrace: \exists f \in F: \langle f, e \rangle \neq 0; ~~~	\forall f \in F \setminus \lbrace 0 \rbrace: \exists e \in E: \langle f, e \rangle \neq 0.
	\end{displaymath}
For instance, a vector space $E$ over $\K$ and its algebraic dual $E^*$ form a dual pair. In fact, fixing $E$, any possible $F$ may be identified with a subset of $E^*$ using the isomorphic embedding $f \rightarrow \langle \cdot, f \rangle$ . Note that the definition of a dual pair is symmetric, so that we can flip the roles of $E$ and $F$ in any definition or proposition. Given a dual pair $(E, F)$, Given $V \subseteq E$, its \textbf{polar} is
\begin{displaymath}
V^{\circ} := \lbrace f \in F: \sup_{e \in V} |\langle e,f \rangle| \leq 1 \rbrace \subseteq F.
\end{displaymath}
The \textbf{bipolar} of $V$ is $V^{\circ\circ} := (V^\circ)^\circ$. It can be proven that for $U, V \subseteq E$:
\begin{itemize}
	\item $U \subseteq V \Rightarrow U^\circ \supseteq V^\circ$ hence $(U \cap V)^\circ \supseteq U^\circ \cup V^\circ$;
	\item $(U \cup V)^\circ = U^\circ \cap V^\circ$;
	\item $V \subseteq V^{\circ\circ}$ hence  $V^{\circ\circ\circ} = V^{\circ}$. 
\end{itemize}

The non-standard extensions $\enl{E}$ and $\enl{F}$ form a dual pair of vector spaces over $\enl{\K}$. Because of the transfer principle $\enl{(V^{\circ})} = (\enl{V})^{\circ}$, hence we may write $\enl{V^{\circ}}$ unambiguously.

\begin{lem} \label{lem_polarKclosed}
	Given a dual pair $(E, F)$ and a set $A \subseteq \enl{E}$ closed for $\stcp\K$-scalar multiplication. Then 
	\begin{displaymath}
	A^{\circ} = \lbrace f \in \enl{F}: |\langle e,f \rangle| \approx 0, \forall e \in A \rbrace.
	\end{displaymath}
	\begin{proof}
		Take $f \in \enl{F}$. Then $f \in A^{\circ}$ iff $|\langle f, \alpha e \rangle| \leq 1$ for any $e \in A$ and $\alpha \in \R_+$, which is equivalent with $|\langle f, e \rangle| \leq \frac{1}{\alpha}$. Since $\alpha$ was arbitrary, $\langle f, e \rangle \approx 0$.
	\end{proof}
\end{lem}

In particular, if $E$ is a TVS then $f \in \m_E(0)^\circ$ iff $e \approx_\tau 0 \Rightarrow \langle e, f \rangle \approx 0$, i.e. $f$ is $\se$continuous (identified as en element of $\enl{F}^*$). If $E$ is a TVS, then $E'$ is the subset of $E^*$ containing the continuous linear functionals, the \textbf{topological dual}. Then we have $\stcp{E'} = \stcp{E}^* \cap \m_E(0)^\circ$.

\begin{lem} \label{lem_polarbalconv}
	Given a dual pair $(E, F)$ and $A \subseteq \enl{F}$. Then $A^\circ$ is $\enl$balanced and $\enl$convex.
	\begin{proof}
		Let $e_1, e_2, \dots, e_\omega$ be an internal sequence in $A^\circ$ and $\alpha_1, \dots, \alpha_\omega $ and internal sequence $\enl{\K}$ such that $\sum_{i=1}^\omega |\alpha_i| \leq 1$ then we get for $f \in A$:
		\begin{displaymath}
		|\langle f, \sum_{i=1}^\omega \alpha_i e_i \rangle| \leq \sum_{i=1}^\omega |\alpha_i| |\langle f, e_i \rangle| \leq 1.
		\end{displaymath}
	\end{proof}
\end{lem}

\begin{lem} \label{lem_polarclosed}
Suppose $E$ is a locally convex Hausdorff TVS. Given a dual pair $(E, F)$ and $A \subseteq \enl{(F \cap {E}')}$. Then $A^\circ$ is $\enl$closed. 
\begin{proof}		
	Take $e$ in the $\enl$closure of $A^\circ$ and $f \in A$. Since $f$ is $\enl$continuous, there is a $\lambda \in \enl\Lambda$ such that $|\langle e', f \rangle| \leq 1$ for any $e' \in \lambda$.  For any $\epsilon \in \enl\R_+$ there exists a $e_\epsilon \in A^\circ$ such that $e - e_\epsilon \in \epsilon\lambda$. Then $|\langle e, f \rangle| \leq |\langle e_\epsilon, f \rangle| + |\langle e - e_\epsilon, f \rangle| \leq 1 + \epsilon$. Since $\epsilon$ was arbitrary, $e \in A^\circ$.
\end{proof}
\end{lem}

\begin{defn}
	Given a dual pair $(E, F)$ where $E$ is a TVS. A set $V \subseteq E$ satisfies the \textbf{bipolar identity} if $V^{\circ\circ} = \overline\cobal(V)$. 
\end{defn}

\begin{lem} \label{lem_dualpairwithbipolar}
	Suppose $E$ is a locally convex Hausdorff TVS. Given $F \subseteq E^*$. Suppose any convex and balanced $V \subseteq E$ satisfies for the duality $(E, E^*)$
	\begin{displaymath}
	(V^\circ \cap F)^\circ = \overline{V}.
	\end{displaymath}
	Then $(E, F)$ is a dual pair for which all sets in $E$ satisfy the bipolar identity.
	
	\begin{proof}
		We know $F$ is separated by $E$. Given $e \in E \setminus \lbrace 0 \rbrace$. There is a $\lambda \in \Lambda$ such that $e \notin \overline\lambda$. The assumption implies $e \notin (\lambda^\circ \cap F)^\circ$, i.e. there exists an $f \in \lambda^\circ \cap F$ such that $\langle e, f \rangle > 1$, in particular $|\langle e, f \rangle| \neq 0$. Now in the duality $(E, F)$, given $W \subseteq E$ lemma \ref{lem_polarbalconv} and order preservation imply that $W^{\circ\circ} = \cobal{(W)}^{\circ\circ}$. But the assumption implies $\cobal{(W)}^{\circ\circ} = \overline\cobal{(W)}$. 
	\end{proof}
\end{lem}

\begin{prop}
	Suppose $E$ is a locally convex Hausdorff TVS. Then $(E, E')$ is a dual pair for which all sets in $E$ satisfy the bipolar identity.
	
	\begin{proof}
		We will apply lemma \ref{lem_dualpairwithbipolar}. Given a convex, balanced $V \subseteq E$. Suppose $e \notin \overline{V}$. Because of regularity, there exists a $\lambda \in \Lambda$ such that $\m(e) \cap (V + \lambda) = \emptyset$, i.e. $e \notin \overline{V + \lambda}$. Since $\overline{V + \lambda}$ is closed, convex, balanced and absorbing we can apply corollary \ref{cor_hbsep}, i.e. there exists a linear functional $f \in E'$ on $E$ such that $|f(e)| > 1$ and $f \in (\overline{V + \lambda})^\circ \subseteq \lambda^\circ \cap V^\circ \subseteq E' \cap V^\circ$. Hence, $e \notin (E' \cap V^\circ)^\circ$. The opposite inclusion follows from lemma \ref{lem_polarclosed}.
	\end{proof}
\end{prop}

By the transferring principle, all internal sets in the duality $(\enl{E}, \enl{E'})$ satisfy the bipolar identity. For external sets, lemma's \ref{lem_polarbalconv} and \ref{lem_polarclosed} apply so that their bipolars are $\enl$closed, $\enl$convex, $\enl$balanced sets; but they can be larger than the hulls. In many cases, the bipolar identity will apply to external sets of interest too, with one notable exception. 

\begin{lem} \label{lem_polgalismon}
	Given a dual pair $(E, F)$, $\mon \in \ps\M({F})$, resp. $\gal \in \ps\Ga({F})$ with $\ps$base $(U_j)_{j \in \enl{J}}$. Then $\mon^\circ \in \ps\Ga({F})$, resp. $\gal^\circ \in \ps\M({F})$ with base $(U_j^\circ)_{j \in \enl{J}}$.
	
	\begin{proof}
		The polar mapping is order-reversing and $\menlp$-extendable.
	\end{proof}
\end{lem} 

\begin{lem} Given a dual pair $(E, F)$. Then \label{lem_bipoolsigma} $(\stcp{E})^{\circ\circ} = \qs(\enl{E})$. 
\begin{proof}		
	Take any $S \in \Pwfin(E)$. Note that $\spann(S)$ is closed for scalar multiplication, so that $\spann(S)^\circ$ and $\spann(S)^{\circ\circ}$ are annihilators. Since $F$ separates $\spann(S)$, we have $	\lbrace \restr{f}{\spann(S)} :  f \in F \rbrace = \spann(S)^*.$
	From this observation, applying well-known principles of linear algebra, it is clear that $\spann(S)^{\circ\circ} = \spann(S)$.
	Let $\lbrace e_1, \dots e_n \rbrace \subseteq S$ be a linearly independent basis of this space. Take $e \in S^{\circ\circ} \subseteq \spann(S)$, i.e. there exists $(\alpha_j)_{j \leq n}$ such that $e = \sum_{j = 1}^n \alpha_j e_j$. Suppose that $|\alpha_k| > 1$ for some $k \leq n$. Applying again separation, there must exist an $f \in F$ such that $\langle e_j, f \rangle = 0$ ($j \neq k$) and $\langle e_k, f \rangle = 1$. Then $f \in S^\circ$ but $|\langle e, f \rangle| = |\alpha_k| > 1$, a contradiction. We conclude that $S^{\circ\circ} = \cobal(S)$. Then the result follows from the fact that the bipolar map is order preserving and thus $\menlp$-extendable.
\end{proof}
\end{lem}

\begin{lem}  \label{lem_mongalbipolar}
Given a dual pair $(E, F)$ with $E$ a TVS. $\ps$Monads and $\ps$galaxies that have a $\ps$base of sets that satisfy the bipolar identity, satisfy the bipolar identity. In particular, for the dual pair $(E, E')$ $\ps$monads and $\ps$galaxies satisfy the bipolar identity.
\begin{proof}
	The bipolar mapping is order preserving and $\menlp$-extendable. If the mapping coincides with $\overline\cobal$ on the sets in the $\ps$base, their $\menlp$-extensions must also coincide on the monad or galaxy.
\end{proof}
\end{lem}

\subsection{Dual topologies}

In this section we work with a duality $(E,F)$ without any explicit topology on $E$.
If $\gal \in \Ga(\enl{F})$ is closed for $\K$-scalar multiplication, then $\gal^\circ$ is closed for $\K_\f$-scalar multiplication and addition by lemma \ref{lem_polarKclosed}. That means we put a topology (also denoted as $\gal$) on $E$ that makes addition continuous, defined by $\m_\gal(0) := \gal^\circ$, i.e. the \textbf{polar topology} of $\gal$ or \textbf{topology of uniform convergence} on $\gal$.\footnote{Recall that continuous addition is sufficient to use the $\m$-equivalence notation from definition \ref{def_mequivalentie}. We will in fact apply this topology to $F^*$ in its entirety ($\m_{F^*}(0) := \gal^\circ$ where we take the polar for the duality $(F^*, F)$), even if it is only compatible with the vector structure on $E$.} This is the coarsest topology that makes every convergent $\mon \in \M(\enl{E})$ uniformly convergent on $\gal$, which is defined as (given $e \in E$):
\begin{displaymath}
\mon \xrightarrow{\gal} e \iff \mon \subseteq e + \gal^\circ \iff \forall V \in \idl \gal: \exists W \in \fil \mon: W - e \subseteq V^\circ
\end{displaymath}
where we applied Cauchy principle on the second equivalence. Because of lemma \ref{lem_polarbalconv} the topology is locally convex. For any base of $\gal$, its poles are a base of closed, balanced, convex neighbourhoods of 0, i.e. of seminorms in the polar topology. We will require additional conditions for this topology to be compatible with the vector structure and Hausdorff. 

\begin{defn}
	The \textbf{weak topology} or topology of pointwise convergence $\sigma(E, F)$ is defined by
	\begin{displaymath}
	\m_{E, \sigma}(0) := (\stcp{F})^\circ.
	\end{displaymath}
\end{defn}

\begin{lem} \label{lem_dualesluitingzwak}
	$F = (E_\sigma)'$.
\begin{proof}		
	We denote $\bullet$ as the polar for the duality $(E, E^*)$. Applying lemma  \ref{lem_bipoolsigma} to this duality, we find:
	\begin{displaymath}
	\stcp{(E_\sigma)'} = \stcp{E^*} \cap \m_{E, \sigma}(0)^\bullet = 
	\stcp{E^*} \cap (\stcp{F})^{\circ\bullet} = 
	\stcp{E^*} \cap \qs(\enl{F}) = \stcp{F}. 
	\end{displaymath}
\end{proof}
\end{lem}

\begin{defn}
	A topology $\tau$ on $E$ is called \textbf{compatible with the duality} if $(E_\tau)' = F$.
\end{defn}

In particular, the weak topology is compatible with the duality. Hence, by lemma \ref{lem_polarclosed} any polar is weakly $\enl$closed. By lemma \ref{lem_mongalbipolar} the bipolar of any $\ps$monad or $\ps$galaxy is equal to the weakly $\enl$closed $\enl$convex $\enl$balanced hull. For topologies compatible with the duality, the closure of any convex balanced set $V \subseteq E$ is equal to the weak closure: $\overline{V}^\tau = \overline{\cobal{(V)}}^\tau = V^{\circ\circ} = \overline{V}^\sigma$.

\begin{prop} \label{prop_taupolair}
	A topology $\tau$ on $E$ is polar for $(E,F)$ iff $\tau$ is locally convex and $\m_\tau(0)$ is weakly $\enl$closed.	
	\begin{proof}
		If $\tau$ is the polar topology of $\gal \in \Ga(\enl{F})$, then $\m_\tau(0) = \gal^\circ$, hence weakly $\enl$closed. If $\m_\tau(0)$ is $\enl$convex and weakly $\enl$closed, then it must be the polar topology of $\m_\tau(0)^\circ$.
	\end{proof}
\end{prop}

\begin{cor}
	For \textit{any} polar topology $\tau$ we have $\m_\tau(0)^{\circ\circ} = \m_\tau(0)$.
\end{cor}

Hence, any polar topology $\tau$ is the topology of uniform convergence on $\m_\tau(0)^\circ$. The polar topology of $\gal$ is also the polar topology of $\overline\cobal(\gal)$, i.e.  uniform convergence on the members of $\idl \gal$ implies uniform convergence on their closed convex balanced hulls. Since $\m_{E, \sigma}(0) = \m_{E, \sigma}(0)^{\circ\circ} =  \qs(\enl{E})^\circ$, the weak topology is the topology of uniform convergence of convex balanced hulls of finite sets.

\begin{prop}
	Given $\gal \in \Ga(\enl{F})$ such that $\overline{\cobal}^\sigma(\gal)$ is covering. The polar topology of $\gal$ is Hausdorff. In particular the weak topology is Hausdorff.
	\begin{proof}
		Assume w.l.o.g. $\gal = \cobal(\gal)$, then $\overline{\gal}^\sigma$ is covering. Take any $e \in \ker \m(0)$, then $\langle e, f \rangle \approx 0$ for each $f \in \gal$; in fact this counts for any $f \in \overline{\gal}^\sigma \supseteq \stcp{F}$ since $e$ is continuous for the weak topology on $F$. But if $f \in \stcp{F}$, $\langle e, f \rangle$ is standard and thus exactly $0$. Hence $e = 0$. 
	\end{proof}
\end{prop}

\begin{prop}	\label{prop_polarfin}
	Given $\gal \in \Ga(\enl{F})$ closed for $\K$-scalar multiplication. Then we have 
	\begin{displaymath}
	\begin{split}
	(\m_\K(0) \cdot \gal)^\circ = \fin_\gal(\enl{E}) \textnormal{~ and ~} \fin_\gal(\enl{E})^\circ &= \m_\K(0) \cdot \overline\cobal(\gal);\\
	\textnormal{in particular }(\m_\K(0) \cdot \stcp{E})^\circ = \fin_\sigma(\enl{E})\textnormal{~ and ~}\fin_\sigma(\enl{E})^\circ &= \m_\K(0) \cdot \qs(\enl{E}).
	\end{split}
	\end{displaymath}
	
	Hence, $\fin_\gal(\enl{E})$ and $\m_\K(0) \cdot \gal$ satisfy the bipolar identity.
	
	\begin{proof}
		Given $e \in \enl{E}$ we find, because of proposition \ref{prop_charfin}:
		\begin{displaymath}	
		\begin{split}
		e \in (\m_\K(0) \cdot \gal)^\circ 
		&\iff \forall \omega \in \enl{\K}_\infty: e \in \omega \cdot \gal^\circ (= \omega \cdot \m_{E, \gal} (0))
		\\&\iff \m_\K(0) \cdot e \subseteq \m_{E, \gal}(0)		
		\\&\iff e \in \fin_\gal(\enl{E}).
		\end{split}
		\end{displaymath}

		Since $\m_\K(0) \cdot \overline\cobal(\gal) \subseteq (\m_\K(0) \cdot \gal)^{\circ\circ}$, by proving the bipolar identity for $\m_\K(0) \cdot \gal$ in case $\gal = \overline\cobal(\gal)$ we get the result. Let $\U$ be a base of $\gal$ existing out of closed, convex, balanced sets and given $\gal' \in \Ga(\enl{F})$. Then $\m_\K(0) \cdot \gal = \bigcup_{U \in \U} \m_\K(0) \cdot \enl{U} \subseteq \gal'$ iff (by the combined Cauchy principle) for any $U \in \U$, there exists an $\alpha \in \K$ and $V \in \idl \gal'$ such that $U \subseteq \alpha V$. In other words, there exist $(\alpha_U)_{U \in \U}$ in $\R_+$ such that $\bigcup_{U \in \U}  \enl{(\alpha_U U)} \subseteq \gal'$, i.e. any galaxy that contains $\m_\K(0) \cdot \gal$, contains a galaxy of this form. Now, since $\m_\K(0) \cdot \gal$ is a union of monads and therefore chromatic, by proposition \ref{prop_chromgal} it is equal to the intersection of all galaxies that contain it. But then applying the previous, 
		\begin{displaymath}
		\m_\K(0) \cdot \gal = \bigcap_{(\alpha_U)_{U \in \U}} \bigcup_{U \in \U} \enl{(\alpha_U U)};
		\end{displaymath}
		an intersection of galaxies that are equal to their bipolar (due to lemma \ref{lem_mongalbipolar}). Then applying lemma \ref{lem_hullintersect} on the $\menlp$-extension of the bipolar map, we find that $\m_\K(0) \cdot \gal$ must be equal to its bipolar.
	\end{proof}
\end{prop}

\begin{cor}
	For any polar topology $\tau$, $\bdd_\tau(\enl{E})^{\circ\circ} = \bdd_\tau(\enl{E})$.
\end{cor}

\begin{cor}
	Given $\mon \subseteq \ps\M(\enl{E})$ and $\gal \in \Ga(\enl{F})$ closed for $\K$-scalar multiplication. Then $\mon$ is $\se$bounded for the polar topology of $\gal$ iff $\gal$ is $\ce$absorbed by $\mon^\circ$. In particular, $\mon$ is weakly $\se$bounded iff $\mon^\circ$ is $\ce$absorbing.
	
	\begin{proof}
		$\mon \subseteq \fin_\gal(\enl{E}) \iff \mon^\circ \supseteq \m_\K(0) \cdot \gal$.
	\end{proof}
\end{cor}

\begin{cor}	
	Given $\gal \in \Ga(\enl{F})$ closed for $\K$-scalar multiplication. Then $\gal^\circ$ is linear iff $\gal$ is weakly $\se$bounded. In particular the weak topology is compatible with the vector structure.
	\begin{proof}
	 	$\gal^\circ \supseteq \m_\K(0) \cdot \stcp{E} \iff \gal \subseteq \fin_\sigma(\enl{E})$.
	\end{proof}
\end{cor}

Hence, the finest polar topology compatible with the vector structure is determined by the weakly bounded sets:
\begin{defn}
	The \textbf{strong topology} $\bt(E, F)$ or topology of uniform convergence on weakly bounded sets is defined by
	\begin{displaymath}
	\m_{E, \bt}(0) := \bdd_\sigma(\enl{E})^\circ.
	\end{displaymath}	
\end{defn}

\begin{prop} \label{prop_polarequicontinu}	
	$\m_E(0)^\circ$ is the galaxy of equicontinuous subsets of $F$. In other words, a set  $H \subseteq F$ is equicontinuous iff $\enl{H} \subseteq \m_E(0)^\circ$ or equivalently $\enl{H}^\circ \supseteq \m_E(0)$.	
	\begin{proof}
		This follows from proposition \ref{prop_equicontinu}.
	\end{proof}
\end{prop}

Hence, closed, convex, balanced hulls of equicontinuous sets are equicontinuous. Each polar topology on $E$ is the topology of uniform convergence on equicontinuous subsets of $F$. 
In fact, each Hausdorff local convex topology is a dual topology for the duality $(E, E')$, i.e. uniform convergence on equicontinuous subsets of $E'$.

\begin{prop} \label{prop_pnsbipool}
	For any polar $\tau$, $\pns_\tau(\enl{E})$ satisfies the bipolar identity.	
	\begin{proof}
		As in the proof of proposition \ref{prop_cobalpns}, we can write $\pns(\enl{E})$ as an intersection of $\enl$closed, $\enl$convex, $\enl$balanced galaxies, so that lemma \ref{lem_hullintersect} implies the result.
	\end{proof}
\end{prop}

\begin{cor}
	For any polar topology $\tau$, $\pcomp_\tau(\enl{E})^{\circ\circ} = \pcomp_\tau(\enl{E})$.
\end{cor}

\begin{prop} For any $\tau$, \label{prop_pnspool} $\pns_\tau(\enl{E})^\circ = \m_{E, \tau}(0)^\circ \cap \stcp{E}^\circ$.
	\begin{proof}		
		Since $\m(0) \cup \stcp{E} \subseteq \pns(\enl{E})$ we have $\pns(\enl{E})^\circ \subseteq \m(0)^\circ \cap \stcp{E}^\circ$, leaving only the opposite implication to prove. Take $f \in \m(0)^\circ \cap \stcp{E}^\circ$ and $e \in \pns(\enl{E})$. Since $f$ is $\se$continuous there exists a $\lambda \in \stcp{\Lambda}$ such that $\langle e', f \rangle \leq 1$ for all $e' \in \lambda$. But there exists a $e_\lambda \in \stcp{E}$ such that $e - e_\lambda \in \lambda$ and therefore
		\begin{displaymath}
		|\langle e, f \rangle| \leq |\langle e
		_{\lambda}, f \rangle| + |\langle e - e_{\lambda}, f \rangle| \lessapprox 1.
		\end{displaymath}
		Since $\pns(\enl{E})$ is closed for $\K$-scalar multiplication $|\langle e, f \rangle| \approx 0$ and therefore $f \in \pns(\enl{E})^\circ$. 	
	\end{proof}
\end{prop}

\begin{cor}
	The weak topology is Henson-Moore, i.e. $\fin_\sigma(\enl{E}) = \pns_\sigma(\enl{E}).$
	
	\begin{proof}
		Take $f \in \pns_\sigma(\enl{E})^\circ = \m_\sigma(0)^\circ \cap \stcp{E}^\circ = \qs(\enl{E}) \cap \m_{F, \sigma}(0)$. Then by the definition of $\qs(\enl{E})$, there are finite $e_0,\dots,e_n \in \stcp{E}$ such that $f = \sum_{j = 0}^n \alpha_j e_j$ with $\alpha_j \approx 0$ for all $j = 0\dots n$. Then for $e \in \fin_\sigma(\enl{E})$, $\langle e, f \rangle = \sum_{j = 1}^n \alpha_j \langle e_j, f \rangle \approx 0$. So $\pns_\sigma(\enl{E})^\circ \subseteq \fin_\sigma(\enl{E})^\circ$ from which follows that $\fin_\sigma(\enl{E}) \subseteq  \pns_\sigma(\enl{E})$ since both sets are equal to their bipolar.
	\end{proof}
\end{cor}

\begin{prop} For any $\tau$, $\ns_\tau(\enl{E})^\circ = \m_{E, \tau}(0)^\circ \cap \stcp{E}^\circ = \pns_\tau(\enl{E})^\circ$.
	\begin{proof}		
		Since $\m(0) \cup \stcp{E} \subseteq \ns(\enl{E})$ we have $\ns(\enl{E})^\circ \subseteq \m(0)^\circ \cap \stcp{E}^\circ$. Then we should only prove the opposite inclusion. Take $f \in \m(0)^\circ \cap \stcp{E}^\circ$ and $e \in \ns(\enl{E})$. There exists a $\hat{e} \in \stcp{E}$ such that $e \approx \hat{e}$ and therefore
		\begin{displaymath}
		|\langle e, f \rangle| \leq |\langle \hat{e}, f \rangle| + |\langle e - \hat{e}, f \rangle| \approx 0.
		\end{displaymath}
	\end{proof}
\end{prop}

\begin{cor}
	For any topology on $E$, on equicontinuous sets of $F$ the topology of uniform convergence on compact sets coincides with the weak topology.
	\begin{proof}
		Indeed, since
		\begin{displaymath}
		\begin{split}
		\m_E(0)^\circ \cap \comp(\enl{E})^\circ &\supseteq \m_E(0)^\circ \cap \ns(\enl{E})^\circ \\&= \m_E(0)^\circ \cap \m_E(0)^\circ \cap \stcp{E}^\circ \\&= \m_E(0)^\circ \cap \m_\sigma(0).		
		\end{split}
		\end{displaymath}
	\end{proof}
\end{cor}

\begin{cor}
	For any polar $\tau$, $\ns_\tau(\enl{E})^{\circ\circ} = \pns_\tau(\enl{E})^{\circ\circ} = \pns_\tau(\enl{E}).$
\end{cor}

In other words, the bipolar identity is only valid for the near-standard points in spaces where completeness and the convex compactness property coincide. Since there exist spaces that have this property but are not complete (for instance quasi-complete spaces that are not complete \cite{Khaleelulla1982}), the bipolar identity is not always valid for the near-standard points. Hence, this is a demonstration of two order-preserving maps that coincide on standard and internal sets, $\enl$monads and $\enl$galaxies but not on all external sets.

\subsection{Compatibility with the duality and completeness}
In this section we work with a duality $(E, F)$ with $\tau$ being a locally convex Hausdorff topology on $E$. Identifying $E$ as a subset of $F^*$, we may extend the definition of the $\st \cdot$ map to $\fin_\sigma(\enl{E}) = \pns_\sigma(\enl{E})$ as follows:
\begin{displaymath}
\langle \stcomt{\sigma} e, f \rangle := \st \langle e, f \rangle \textnormal{ for } f \in \stcp{F}.
\end{displaymath}
Then $\stcomt{\tau} \cdot$ is a restriction of this map to $\pns_\tau(\enl{E})$. Indeed, since $\stcomt{\sigma} e \approx_\sigma e$ and the topology of $E$ is finer than the weak topology, this extension is entirely compatible with the previous definition of $\stcom$. 
\begin{lem} \label{lem_monadsatdual}
	Given $V \subseteq E$ a convex balanced set. Then
	\begin{displaymath}
	(\m_{E, \tau}(0) \cap \enl{V})^\circ \subseteq \stcp{({V_\tau}' \cap E^*}) + \m(0) = \ns_\sigma(\enl{({V_\tau}' \cap E^*)}).
	\end{displaymath}
	In other words, any member of $\enl{F}$ that is $\se$continuous on $V$ is in the weak topological monad (for the space $E^*$) of a functional that is continuous on $V$.
	\begin{proof}
		Take $f \in (\enl{V} \cap \m_E(0))^\circ$. Because of the Cauchy principle, there exists a $U \in \fil \m_E(0)$ so that $f \in \enl(V \cap U)^\circ$. Let $\hat{f} := \stcomt{\sigma} f \in E^*$. Then we have for any $e \in U$ that $|\langle e, \hat{f} \rangle| \approx |\langle e, f \rangle| \leq 1$ so that $|\langle e, \hat{f} \rangle| \leq 1$ since it has a standard value, i.e. $\hat{f} \in (V \cap U)^\circ \subseteq V'$. 
	\end{proof}
\end{lem}

\begin{lem} \label{lem_monadsatdual2}
	Given $f \in \stcp{(V_\tau)' \cap E^*}$ with $V \subseteq E$ a convex balanced set, there exists an $f' \in (\m_{E, \tau}(0) \cap \enl{V_\tau})^\circ$ such that that $f' \approx_\sigma f$.
	In other words, any functional in $E^*$ that is continuous for $V$, is in the weak topological monad of a member of $\enl{F}$ that is $\se$continuous on $V$.
	\begin{proof}
		Notice that $f \in (\m_E(0) \cap \enl{V_\tau})^\bullet$, where $\bullet$ is the polar for the duality $(E, E^*)$. Since $\m_E(0)^\bullet = \m_E(0)^{\circ\circ\bullet} = \m_E(0)^{\circ\bullet\bullet}$, this is in fact the weak closure in $\enl{E^*}$ (for the duality $(E, E^*)$) of $\m_E(0)^\circ$, implying the result.
	\end{proof}
\end{lem}

\begin{thm} \label{th_monadpolarns}
 $\tau$ is compatible with the duality iff $\tau \geq \sigma$ and
	\begin{displaymath}
	\m_{E, \tau}(0)^\circ \subseteq \ns_\sigma(\enl{F}).
	\end{displaymath}
	
	\begin{proof}
		\begin{displaymath}
		F \subseteq E' \iff \stcp{F} \subseteq \m_E(0)^\circ \iff \m_E(0) \subseteq \stcp{F}^\circ = \m_{E,\sigma}(0).
		\end{displaymath}			
		We must then only proof that $\m_{E, \tau}(0)^\circ \subseteq \ns_\sigma(\enl{F})$ is equivalent with $E' \subseteq F$. The $\Leftarrow$-implication is lemma \ref{lem_monadsatdual} (with $V = E$). The $\Rightarrow$-implication follows from lemma \ref{lem_monadsatdual2}. Indeed, given $f \in E'$, we find a $f' \in \m_{E, \tau}(0)^\circ$ so that $f \approx_\sigma f$ which implies that $\st f' = f$ and therefore $f \in F$.
	\end{proof}
\end{thm}

\begin{cor} \label{cor_albou}
	(Alaoglu-Bourbaki) The equicontinuous sets of $E'$ are weakly relatively compact.
\end{cor}

\begin{cor} (Mackey-Arens)
	A locally convex Hausdorff topology $\tau$ on $E$ is compatible with the duality iff $\tau$ is a topology of uniform convergence on covering of convex, balanced, weakly relatively compact sets. 
\end{cor}

The theorem is thus the non-standard characterization of both Alaoglu-Bourbaki and Mackey-Arens which now appear to be trivially equivalent. The most coarse topology that is compatible with the duality is the weak topology. Now we know from section \ref{sec_coco} that there exists a greatest galaxy and therefore a smallest monad or finest topology $\tau(E, F)$ that satisfies this condition, namely $\m_{E, \tau}(0) = \coco_\sigma(\enl{F})^\circ$, the \textbf{Mackey-topology}.

\begin{thm}
	$\pns_\tau(\enl{E}) = (\m_{E, \tau}(0)^\circ \cap \stcp{E}^\circ)^\circ$.
	
	\begin{proof}
		This follows from propositions \ref{prop_pnsbipool} and \ref{prop_pnspool}.
	\end{proof}
\end{thm}

\begin{cor} \label{cor_pnscontinu}
	Given $\U$ a subbase of $\m_{E, \tau}(0)^\circ$. Then
	$e \in \pns_\tau(\enl{E})$ iff $e$ is weakly $\se$continuous on each $U \in \U$ (weakly referring to $\sigma(F, E)$).
	
\begin{proof}
	$\pns(\enl{E}) = (\m_E(0)^\circ \cap \stcp{E}^\circ)^\circ = (\m_E(0)^\circ \cap \m_{F, \sigma}(0))^\circ$
	\begin{displaymath}
	= (\bigcup_{U \in \U} \enl{U} \cap \m_{F, \sigma}(0) )^\circ
	= \bigcap_{U \in \U} (\enl{U} \cap \m_{F, \sigma}(0) )^\circ.
	\end{displaymath}
\end{proof}	

\end{cor}

\begin{cor} (Grothendieck completeness theorem) Given $\U$ a subbase of $\m_E(0)^\circ$. 
	Then $E$ is complete iff each $e \in F^*$ that is weakly continuous on elements of $\U$ is weakly continuous on the entire space. The completion of $E$ is precisely
	\begin{displaymath}
	\hat{E} = \lbrace e \in F^*: \forall U \in \U: \restr{e}{U}\textnormal{ is weakly continuous} \rbrace
	\end{displaymath}	
	\begin{proof}
		The characterization of $\hat{E} = \stcomt{\sigma} \pns(\enl{E})$ follows from corollary \ref{cor_pnscontinu} combined with lemmas \ref{lem_monadsatdual} and \ref{lem_monadsatdual2}. Lemma \ref{lem_dualesluitingzwak} implies $E = \stcomt{\sigma} \ns(\enl{E}) = (F_\sigma)'$.
	\end{proof}
	
\end{cor}

\bibliographystyle{plainnat} 
\bibliography{dual.bib}

\end{document}